\pdfoutput=1
\RequirePackage{ifpdf}
\ifpdf 
\documentclass[pdftex]{sigma}
\else
\documentclass{sigma}
\fi

\numberwithin{equation}{section}

\usepackage[all,cmtip]{xy}

\RequirePackage{tikz-cd}
\usetikzlibrary{calc}
\usetikzlibrary{decorations.pathmorphing}

\tikzset{curve/.style={settings={#1},to path={(\tikztostart)
			.. controls ($(\tikztostart)!\pv{pos}!(\tikztotarget)!\pv{height}!270:(\tikztotarget)$)
			and ($(\tikztostart)!1-\pv{pos}!(\tikztotarget)!\pv{height}!270:(\tikztotarget)$)
			.. (\tikztotarget)\tikztonodes}},
	settings/.code={\tikzset{quiver/.cd,#1}
		\def\pv##1{\pgfkeysvalueof{/tikz/quiver/##1}}},
	quiver/.cd,pos/.initial=0.35,height/.initial=0}

\tikzset{tail reversed/.code={\pgfsetarrowsstart{tikzcd to}}}
\tikzset{2tail/.code={\pgfsetarrowsstart{Implies[reversed]}}}
\tikzset{2tail reversed/.code={\pgfsetarrowsstart{Implies}}}
\tikzset{no body/.style={/tikz/dash pattern=on 0 off 1mm}}

\newtheorem{Thm}{Theorem}[section]
\newtheorem{Lem}[Thm]{Lemma}
\newtheorem{Con}[Thm]{Condition}
\newtheorem{Prop}[Thm]{Proposition}
\newtheorem{Cor}[Thm]{Corollary}

\theoremstyle{definition}
\newtheorem{Def}[Thm]{Definition}
\newtheorem{Exm}[Thm]{Example}
\newtheorem{Rem}[Thm]{Remark}

\begin{document}

\allowdisplaybreaks

\newcommand{\arXivNumber}{2412.12935}

\renewcommand{\PaperNumber}{074}

\FirstPageHeading

\ShortArticleName{$\mathcal{L}$-Lie Algebroids over Topological Ringed Spaces}

\ArticleName{$\boldsymbol{\mathcal{L}}$-Lie Algebroids over Topological Ringed Spaces}

\Author{Mainak PODDAR~$^{\rm a}$ and Abhishek SARKAR~$^{\rm ab}$}

\AuthorNameForHeading{M.~Poddar and A.~Sarkar}

\Address{$^{\rm a)}$~Department of Mathematics, Indian Institute of Science Education and Research Pune, India}
\EmailD{\mail{mainak@iiserpune.ac.in}}

\Address{$^{\rm b)}$~School of Computing, MIT Art, Design and Technology University, Loni Kalbhor, Pune, India}
\EmailD{\mail{abhishek.sarkar@mituniversity.edu.in}}

\ArticleDates{Received October 14, 2025, in final form August 01, 2026; Published online August 11, 2026}

\Abstract{Lie algebroids over a topological ringed space provides a unified framework to study various geometric structures. This geometric concept is intimately connected with well-known algebraic structures, including Gerstenhaber algebras and Batalin--Vilkovisky algebras. We introduce more general concepts such as $\mathcal{L}$-Lie algebroids and $\mathcal{A}$-Gerstenhaber algebras, associated with a given Lie algebroid $\mathcal{L}$ and Gerstenhaber algebra $\mathcal{A}$ over a topological ringed space, respectively. Following this, we explore how several standard correspondences extend within this broader framework.}

\Keywords{Lie algebroids; tangent sheaf; Gerstenhaber algebras; BV-algebras; universal enveloping algebra; Lie algebroid connections; Lie algebroid (co)homology}

\Classification{53D17; 32L10; 14F40; 17B66}

\section{Introduction}

Lie algebroids plays a prominent role in differential geometry and mathematical physics as they represent generalized infinitesimal symmetries of spaces, which are related to the corresponding global symmetries of spaces described by Lie groupoids~\cite{KM}.
Lie algebroids over a~$C^{\infty}$\nobreakdash-manifold generalize the tangent bundle of a manifold as well as Lie algebras \cite{AW-EL, MM}.
The standard geometric notions (or structures) associated with a Lie algebroid such as exterior algebra of a Lie algebroid, Lie bialgebroid, Lie algebroid connections, Chevalley--Eilenberg--de Rham complex, etc., are defined in the context of smooth vector bundles (see~\cite{KM, PX}). In recent years, some of these notions have been generalized in the context of complex geometry and algebraic geometry, using sheaf theoretic language (see \cite{UB,BRT,DRV, CV,MK, AA, BP, PT}). These geometric notions are closely associated with some well-known algebraic structures, such as Lie--Rinehart algebras, Gerstenhaber algebras (in short G-algebras), Batalin--Vilkovisky algebras (in short BV-algebras), universal enveloping algebra of a Lie--Rinehart algebra, etc.\ (see \cite{JH, JH2,KS, MM, LM, PX}).
We consider sheaf theoretic (or algebro-geometric) analogues of these algebraic notions, and investigate some of the important relationships between them (see Section~\ref{Sec 2}).\looseness=1

Unlike smooth vector bundles over a smooth manifold, the global sections of a holomorphic or algebraic vector bundle over a complex manifold or algebraic variety respectively, do not, in general, capture all the information about the bundle.
Consequently, a sheaf-theoretic approach becomes essential in analytic and algebraic geometry, even in non-singular cases. To create a~unified framework for studying geometric objects across diverse contexts, such as smooth, holomorphic, and algebraic settings, one relies on sheaves over locally ringed spaces (see~\cite{TS, CW}).
When a space has singularities, the tangent bundle is replaced by the tangent sheaf, which consists of derivations acting on the structure sheaf of the space. Another important sheaf associated to these spaces is the logarithmic tangent sheaf, meaning the sheaf of logarithmic derivations corresponding to some divisor. Both of these notions share a common underlying structure as (quasi)coherent sheaves of Lie--Rinehart algebras, generally regarded as Lie algebroids (see \cite{AA, BP, JV}).

Recently, $L$-Lie algebroids have found several applications in the literature. This structure was introduced in the smooth setting by R.L.~Klaasse as Lie algebroids whose anchors factor through the anchor of a fixed Lie algebroid~$L$~\cite{RK}.
 He has used it in the study of log symplectic structures taking $L$ to be a logarithmic tangent sheaf. An analogous notion was introduced in~\cite{RL25} to extend the concept of singular subalgebroids of a Lie algebroid~\cite{MZ}.
The logarithmic Atiyah sheaf and the $G$-Atiyah sheaf employed by Biswas et al.\ \cite{IBMP, IB, BKN} in the context of equivariant geometry are also examples of $L$-Lie algebroids. Note that when $L$ is the tangential Lie algebroid, an $L$-Lie algebroid is simply a Lie algebroid.

An $L$-Lie algebroid may be informally thought of as a rooted Lie algebroid. In this article, we formulate analogous rooted versions of the notions of differential graded algebras, Gerstenhaber algebras, strong differential Gerstenhaber algebras,
BV-algebras, Lie bialgebroids, etc. We always work in the setting of sheaves over locally ringed spaces. We generalize the relationships among the classical versions of these objects to their rooted avatars. Some of these are given a~functorial interpretation.

Our basic object of study is an $\mathcal{L}$-Lie algebroid, where $\mathcal{L}$ refers to a given Lie algebroid over a~topological ringed space (see Section~\ref{Sec 3}). This notion appears in various geometric settings, such as Poisson geometry \cite{BP}, log symplectic structures \cite{RK}, singular foliations~\cite{CL}, and triangular Lie bialgebroids \cite{DC, KS}.
Recall that, given a Lie algebroid $\mathcal{L}$, there is a sheaf of differential graded algebras $\Omega^\bullet _\mathcal{L}$, known as the Chevalley--Eilenberg--de Rham complex of $\mathcal{L}$ \cite{UB, AS}. Moreover, there is a sheaf of Gerstenhaber algebras $\mathcal{G}_{\mathcal{L}}$, given by the canonical Gerstenhaber algebra structure on the exterior algebra $\wedge^\bullet_{\mathcal{O}_X} \mathcal{L}$ endowed with the generalized Schouten--Nijenhuis bracket \cite{CV, TS}.
The smooth and algebraic counterparts of the above associations are described in \cite{JH1, KS-M, PX}. We define more general notions, namely, $\Omega^\bullet _\mathcal{L}$-differential graded algebras and $\mathcal{G}_{\mathcal{L}}$-Gerstenhaber algebras for a given Lie algebroid $\mathcal{L}$ over a topological ringed space $(X, \mathcal{O}_X)$. We establish analogous relationships among them (see Theorems~\ref{L-version of G-alg & Lie alg} and~\ref{L-Lie alg & dga}).

There is a well-known correspondence between Lie bialgebroid structures on a smooth vector bundle and strong differential Gerstenhaber algebra structures on its associated exterior algebra, due to Kosmann-Schwarzbach (see \cite{JH1,KS, PX}). A~strong differential Gerstenhaber algebra is a~Gerstenhaber algebra together with a square zero differential of degree $1$ which acts like a graded derivation with respect to both the exterior product and bracket operations.
We introduce $\mathcal{L}$-Lie bialgebroids and associate them with strong differential $\mathcal{G}_{\mathcal{L}}$-Gerstenhaber algebras
(see Theorem~\ref{L-version of Lie bialg & strong G alg}). An $\mathcal{L}$-Lie bialgebroid is a pair of $\mathcal{L}$-Lie algebroids in duality with an additional compatibility condition as in the usual case.

For an $\mathcal{O}_X$-module $\mathcal{E}$, the Atiyah algebroid $\mathcal{A}{\rm t}(\mathcal{E})$ consists of first order differential operators on $\mathcal{E}$ having scalar symbol \cite{UB, MK}.
An $\mathcal{L}$-connection $\nabla$ on an $\mathcal{O}_X$-module $\mathcal{E}$ is an $\mathcal{O}_X$-linear morphism $\nabla\colon \mathcal{L} \rightarrow \mathcal{A}{\rm t}(\mathcal{E})$ which satisfies the usual Leibniz rule.
We extend $\mathcal{L}$-connections and their curvature on an $\mathcal{O}_X$-module $\mathcal{E}$ to $(\mathcal{L}', \mathcal{L})$-connections and their curvature, for an $\mathcal{L}$-Lie algebroid $\mathcal{L}'$.
For a smooth Lie algebroid $L$ with a BV-algebra structure on the space of sections of the exterior bundle $\wedge^\bullet L$, P.~Xu defined Lie algebroid homology with coefficients in the line bundle $\wedge^n L$, where $\operatorname{rank} L=n$ (see~\cite{PX}).
We consider a modified version of BV-algebra structure on $\wedge^\bullet_{\mathcal{O}_X} \mathcal{L}$ associated with a $G$-algebra structure on $\wedge_{\mathcal{O}_X}^{\bullet} \mathcal{L}'$ (see Theorem~\ref{L-BV alg}). We associate it with an $\mathcal{L}$-Lie algebroid structure on $\mathcal{L}'$ with a flat $(\mathcal{L}', \mathcal{L})$-connection
on $\wedge^n_{\mathcal{O}_X} \mathcal{L}$, where $\operatorname{rank}(\mathcal{L})=n$. This yields a homology theory for $\mathcal{L}$-Lie algebroids.

As a consequence of the BV-algebra structure on the exterior algebra of a smooth Lie algebroid, a duality theorem between its homology and cohomology was initially established by P.~Xu~\cite{PX} for a specific case and generalized by J.~Huebschmann \cite{JH1, JH2} to the broader framework of Lie--Rinehart algebras. Building on these advancements, we extend this duality correspondence (see Theorem~\ref{L-version duality}), from the context of finitely generated projective Lie--Rinehart algebras to the context of $\mathcal{L}$-Lie algebroids.

The universal enveloping algebroid of a Lie algebroid serves as an analogue of the sheaf of differential operators on a manifold (see \cite{DRV, AA, SR, TS}). It is equipped with a canonical filtered $\mathbb{K}_X$-algebra structure as well as an $\mathcal{O}_X$-bimodule structure (see \cite{DRV, CV, AA}).
This notion plays a pivotal role in exploring the homological algebra and representation theory of Lie algebroids (see \cite{UB, BRT, BP, AS}). A key ingredient in this is a version of the Poincar\'e--Birkhoff--Witt (PBW) theorem for Lie--Rinehart algebras, first established by G.~Rinehart \cite{GR}. Additionally, a sheaf-theoretic analogue of this result for Lie algebroids has been developed (see \cite{MK, AA, BP, TS}).
We extend this framework to incorporate $\mathcal{L}$-Lie algebroids. Specifically, we generalize the PBW theorem of Lie algebroids to the setting of $\mathcal{L}$-Lie algebroids (see Theorem~\ref{L-version PBW}).

In the subsequent sections of this article, we present several examples originating from singular foliations \cite{RF, LGLRLR, CL, BP}, Poisson analytic spaces \cite{BP}, Atiyah algebroids \cite{UB, MK}, free Lie algebroids \cite{MK}, logarithmic connections \cite{IB, DKP}, logarithmic de Rham complex \cite{CMD, AS}, Poisson cohomology \cite{JH, BP}, generalized holomorphic bundles \cite{MG-GC, NH}, equivariant principal bundles \cite{IBMP, IB, BK-MP}, and related structures. While these can be viewed as instances of Lie algebroids over topological ringed spaces, the general $\mathcal{L}$-Lie algebroid framework often provides a simpler and more natural description. Depending on the context, it is sometimes more appropriate to replace the tangential Lie algebroid with another fixed Lie algebroid~$\mathcal{L}$ that plays the central role.

Here we describe a guiding example coming from log-Poisson geometry, which appears consistently throughout the article.
Let $\mathcal{L}:=\mathcal{T}_X(-\log Y)$ be the log tangent sheaf along an analytic hypersurface $Y$ of a holomorphic Poisson manifold $(X,\pi)$. Then the sheaf of meromorphic differential $1$-forms $\mathcal{L}':=\Omega_X^1(\log Y)$ of $X$ with poles of order at most $1$ along $Y$ forms an $\mathcal{L}$-Lie algebroid, and consequently $(\mathcal{L}',\mathcal{L})$ forms an $\mathcal{L}$-Lie bialgebroid over $(X,\mathcal{O}_X)$ (see Example~\ref{log Poisson}).
Here~$\mathcal{G}_{\mathcal{L}}$ is the sheaf of multivector fields on $X$ tangent to $Y$ (or, sheaf of multiderivations of~$\mathcal{O}_X$ preserving the ideal sheaf associated to $Y$) equipped with the canonical sheaf of G-algebras structure (see Example~\ref{multivector fields}). In this situation, we obtain a canonical strong differential G-algebra structure over $(X,\mathcal{O}_X)$ associated with the log Poisson complex $\mathcal{T}^{\bullet}_X(-\log Y)$ of~$X$ along~$Y$ (see Example~\ref{log Poisson strong diff G-algs}).
Moreover, any logarithmic connection on an $\mathcal{O}_X$-module $\mathcal{E}$ naturally produces an $(\mathcal{L}',\mathcal{L})$-connection on $\mathcal{E}$. Additionally, the connection is flat if and only if there exists a~$\mathcal{A}{\rm t}(\mathcal{E})(-\log Y)$-Lie algebroid structure on $\mathcal{L}$ (see Example~\ref{log connection}). Finally, the universal enveloping algebroid $\mathcal{U}(\mathcal{O}_X,\mathcal{L})$ is the sheaf of logarithmic differential operators of $\mathcal{O}_X$ preserving the ideal sheaf associated to $Y$ (see Example~\ref{sheaf of log diff operators}).

\section[Lie algebroids over (X,O\_X)]{Lie algebroids over $\boldsymbol{(X, \mathcal{O}_X)}$} \label{Sec 2}

In this section, we recall some standard notions, namely topological ringed spaces, Lie--Rinehart algebras, Lie algebroids, Gerstenhaber algebras, Chevalley--Eilenberg--de Rham complex, Lie bialgebroids, BV-algebras and so on, in the algebro-geometric language. Then, we recall some important relationships among them in this setting. In the later sections, we extend these notions by considering a rooted version and the analogous relationship in a more general setting.

\subsection{Lie algebroids over topological ringed spaces}
To define topological ringed spaces and Lie algebroids over them, we establish the following notation. Let $\mathbb{K}$ denote either $\mathbb{R}$ or $\mathbb{C}$ (the real or complex number fields, respectively) or, an algebraically closed field of characteristic zero. The constant sheaf on a base space $X$, with stalks isomorphic to $\mathbb{K}$, is denoted by $\mathbb{K}_X$. The sheaf of $\mathbb{K}$-valued continuous functions on $X$, where $\mathbb{K}$ is equipped with either the Euclidean or Zariski topology, is represented by $C_X^0$.

\subsubsection{Topological ringed spaces}
Let $X$ be a topological space and $\mathcal{O}_X$ be a $\mathbb{K}_X$-subalgebra of the sheaf of $\mathbb{K}$-algebras $C_X^0$.
Then the pair $(X, \mathcal{O}_X)$ is said to be a \emph{topological ringed space}.

\begin{Rem} \label{special spaces}
	An appropriate choice for the sheaf $\mathcal{O}_X$ (known as the structure sheaf) over $X$ provides various examples of topological ringed spaces $(X, \mathcal{O}_X)$ coming from smooth manifolds, complex manifolds, analytic spaces and algebraic varieties. These topological ringed spaces are special cases of locally ringed spaces \cite{CW}, we refer to such spaces as special spaces.
\end{Rem}
\begin{Rem} \label{Vec bun & loc free}
	Let $(X, \mathcal{O}_X)$ be a special ringed space. There is a one-to-one correspondence between vector bundles over $X$ (in the respective category) and locally free sheaves of $\mathcal{O}_X$\nobreakdash-modules of finite rank. In fact, the correspondence is an equivalence of categories (see \cite{BP, SR, CW}).
\end{Rem}	

For a topological ringed space $(X, \mathcal{O}_X)$, the sheaf of derivations of $\mathcal{O}_X$-modules is denoted by $\mathcal{D}{\rm er}_{\mathbb{K}_X}(\mathcal{O}_X)$. This is the subsheaf (of $\mathcal{E}{\rm nd}_{\mathbb{K}_X}(\mathcal{O}_X)$) consisting of those endomorphisms on $\mathcal{O}_X$ which satisfy the Leibniz rule.
More precisely, for an open set $U$ of $X$,
\[
\mathcal{D}{\rm er}_{\mathbb{K}_X}(\mathcal{O}_X)(U):= \{D\colon \mathcal{O}_U \rightarrow \mathcal{O}_U \mbox{ is a derivation}\} \subset \mathcal{E}{\rm nd}_{\mathbb{K}_X}(\mathcal{O}_U),
\]
 where $\mathcal{O}_U=\mathcal{O}_X|_U$. That is, a section $D$ of $\mathcal{D}{\rm er}_{\mathbb{K}_X}(\mathcal{O}_X)$ on $U$ is a $\mathbb{K}$-linear derivation~$D$ of~$\mathcal{O}_X(U)$ (or, $D \in \mathcal{D}{\rm er}_{\mathbb{K}}(\mathcal{O}_X(U)))$ such that
for open subsets $W$, $V$ of $U$ with $W \subset V$, satisfies the compatibility condition:
\[
\operatorname{res}_{VW}^{\mathcal{O}_X} \circ D_V= D_W \circ \operatorname{res}_{VW}^{\mathcal{O}_X},
\]
 where $D_{V'}= D|_{V'}$, i.e., \smash{$D_{V'}=\operatorname{res}_{UV'}^{\mathcal{E}{\rm nd}_{\mathbb{K}_X}(\mathcal{O}_X)} (D)$}, for an open subset $V'$ of $U$.

The sheaf $\mathcal{D}{\rm er}_{\mathbb{K}_X}(\mathcal{O}_X)$ over such a space $(X, \mathcal{O}_X)$ has a canonical $ \mathcal{O}_X$-module and $\mathbb{K}_X$-Lie algebra structure, which are compatible via a Leibniz rule \cite{MK,BP, SR}. For the special spaces, it is known as the tangent sheaf of~$X$, and denoted by $\mathcal{T}_X$.

\begin{Rem} \label{nonsingular}
	A special topological ringed space $(X,\mathcal{O}_X)$ is said to be smooth or nonsingular if its tangent sheaf $\mathcal{D}{\rm er}_{\mathbb{K}_X}(\mathcal{O}_X)$ is a locally free $\mathcal{O}_X$-module of finite rank.
\end{Rem}

\subsubsection{Lie algebroids}
Classically, Lie algebroids are defined on smooth (real) vector bundles over a smooth mani\-fold~\cite{KM}. The space of global sections of a Lie algebroid forms a Lie--Rinehart algebra, capturing the entire structure~\cite{JH}.
In contrast, for holomorphic (or algebraic) vector bundles over complex manifolds (or algebraic varieties), the space of global sections does not capture the whole information.

Recall that, a $(\mathbb{K}, A)$-\emph{Lie--Rinehart algebra} $L$ is a $\mathbb{K}$-Lie algebra $L$ with an $A$-module structure (where $A$ is a commutative unital $\mathbb{K}$-algebra) and an $A$-module homomorphism $\rho\colon L\rightarrow \mathcal{D}{\rm er}_{\mathbb{K}}(A)$ such that
\begin{itemize}\itemsep=0pt
	\item the map $\rho$ is also a $\mathbb{K}$-Lie algebra homomorphism,
	\item $ [x, ry]=r[x,y]+\rho(x)(r)y$ for all $x,y\in L$, $r\in A$.
\end{itemize}

\begin{Rem}
	For a Lie algebroid $L$ over a smooth manifold $M$, the space of sections $\Gamma(L|_U)$ of $L$ on an open set $U \subset X$, forms an $(\mathbb{R}, C^{\infty}(U))$-Lie--Rinehart algebra. In addition, the underlying $C^{\infty}(U)$-module structure on $\Gamma(L|_U)$ is finitely generated $($f.g.$)$ projective. Moreover, a f.g.\ projective $(\mathbb{R}, C^{\infty}(M))$-Lie--Rinehart algebra produces a smooth Lie algebroid over $M$.
(Using an analogue of the Serre--Swan theorem, see~\cite{JN}.)
\end{Rem}
A sheaf-theoretic framework is essential in holomorphic and algebraic geometry, as it applies uniformly to both singular and non-singular spaces.
In this setting, we recall Lie algebroids over topological ringed spaces, which generalize both the tangent sheaf and sheaves of Lie algebras (see~\cite{AA, JV}).

\begin{Def}[Lie algebroids over a topological ringed space] \label{Lie algebroids}
	A Lie algebroid $\mathcal{L}$ over a~topological ringed space $(X,\mathcal{O}_X)$ is a sheaf of $\mathbb{K}$-Lie algebras on $X$ and an $\mathcal{O}_X$-module equipped with a homomorphism $\mathfrak{a}\colon (\mathcal{L},[\cdot,\cdot])\rightarrow (\mathcal{D}{\rm er}_{\mathbb{K}_X}(\mathcal{O}_X),[\cdot,\cdot]_c)$ of $\mathcal{O}_X$-modules and $\mathbb{K}_X$-Lie algebras, called the anchor map. The map $\mathfrak{a}$ satisfies the Leibniz rule:
\[
[D,f D'] = f[D,D']+\mathfrak{a}(D)(f) D'\qquad \text{for all} \ f \in \mathcal{O}_X \ \text{and}
\ D,D' \in \mathcal{L}.
\]
\end{Def}	

Therefore, a Lie algebroid $\mathcal{L}$ over $(X,\mathcal{O}_X)$ can be viewed as a sheaf of $(\mathbb{K}_X, \mathcal{O}_X)$-Lie--Rinehart algebras.	We denote a Lie algebroid by $(\mathcal{L}, [\cdot, \cdot], \mathfrak{a})$ or simply by $\mathcal{L}$.

A morphism of Lie algebroids over a topological ringed space
\[
\phi \colon \ (\mathcal{L}_1,[\cdot,\cdot]_1,\mathfrak{a}_1)\rightarrow (\mathcal{L}_2,[\cdot,\cdot]_2,\mathfrak{a}_2)
\]
 is a sheaf homomorphism of Lie--Rinehart algebras, i.e.,
\begin{itemize}\itemsep=0pt
	\item $\phi\colon \mathcal{L}_1 \rightarrow \mathcal{L}_2$ is an $\mathcal{O}_X$-linear map,
	\item $\phi\colon (\mathcal{L}_1,[\cdot,\cdot]_1)\rightarrow (\mathcal{L}_2,[\cdot,\cdot]_2)$ is a $\mathbb{K}_X$-Lie algebra homomorphism,
	\item compatibility: $\mathfrak{a}_2 \circ \phi= \mathfrak{a}_1$.
\end{itemize}

\begin{Rem} \label{loc free f.r.}
	For a locally free $\mathcal{O}_X$-module $\mathcal{E}$ of finite rank $n$, the following properties hold:%
	\begin{enumerate}\itemsep=0pt
		\item For each $x \in X$, there exists an open neighborhood $U_x$ of $x$ in $X$ such that \smash{$\mathcal{E}|_{U_x} \cong \mathcal{O}_X|_{U_x}^{\oplus n}$}.
		\item For each $x \in X$, the stalk $\mathcal{E}_x$ is a free $\mathcal{O}_{X,x}$-module of rank $n$.
	\end{enumerate}
	Moreover, if $\mathcal{E}$ is equipped with a Lie algebroid structure, this induces a Lie--Rinehart algebra structure on the space of sections of $\mathcal{E}$ as well as on its stalks, retaining the above properties.\looseness=1
\end{Rem}
Using a standard approach, for a (locally free) Lie algebroid $\mathcal{L}$ over $(X, \mathcal{O}_X)$, we mainly apply the algebraic properties of the $(\mathbb{K}, \mathcal{O}_X(U))$-Lie--Rinehart algebra $\mathcal{L}(U)$ (which are free or projective module over $\mathcal{O}_X(U)$, for certain open subsets $U \subset X)$ at the local level to derive global results (see \cite{UB, BRT, DRV, CV, MK, AA}).
\begin{Exm}
	A Lie algebroid over a smooth manifold (commonly referred to as a smooth Lie algebroid) $X$ is equivalent to a locally free Lie algebroid of finite rank over the topological ringed space $(X, C_X^\infty)$, where $C_X^\infty$ denotes the sheaf of real valued smooth functions on $X$ (see Remark~\ref{Vec bun & loc free}).
	
A Lie algebroid over a complex manifold $X$ that is locally free of finite rank as an $\mathcal{O}_X$-module ($\mathcal{O}_X$ is the sheaf of holomorphic functions on $X$) is the equivalent to a holomorphic Lie algebroid (see \cite{PT2, BP}).
\end{Exm}

\begin{Exm} \label{Example8}
The standard Lie algebroid structure on the tangent sheaf of a smooth manifold (or complex manifold) $(X, \mathcal{O}_X)$ is described by the triple $(\mathcal{D}{\rm er}_{\mathbb{K}_X}(\mathcal{O}_X), [\cdot,\cdot]_c, {\rm id})=:\mathcal{T}_X$, where $\mathcal{O}_X$ is the sheaf of smooth $($or holomorphic$)$ functions on $X$. The sheaf of sections ${\mathfrak{X}_X:=\Gamma_X(TX)}$ of the smooth $($or holomorphic$)$ tangent bundle~$TX$, i.e., the sheaf of smooth $($or holomorphic$)$ vector fields on~$X$, is isomorphic to $\mathcal{D}{\rm er}_{\mathbb{K}_X}(\mathcal{O}_X)$, forming a locally free Lie algebroid of rank equal to the $($real or complex$)$ dimension of the manifold~$X$~\cite{SR}.
\end{Exm}

\begin{Exm} \label{Foliation}
A (possibly singular) foliation $\mathcal{F}$ on a real smooth manifold or a complex manifold $(X, \mathcal{O}_X)$ is an $\mathcal{O}_X$-submodule of the Lie algebroid $\mathcal{T}_X$ $($or $\Gamma_X(TX))$, which is $(i)$ stable under the Lie bracket and $(ii)$ locally finitely generated. It provides a Stefan--Sussmann distribution or Nagano's distribution on $X$ respectively, and forms a Lie algebroid over $(X, \mathcal{O}_X)$.
	These are coherent involutive subsheaves of the tangent sheaf.
	In particular, regular foliations arise from involutive subbundles of the tangent bundle $($providing Frobenius distributions$)$ are such examples. See \cite{IAGS, RF, LGLRLR, CL, BP} for more details.
\end{Exm}
\begin{Exm} \label{Cotantgent Poisson}
	The sheaf of differential $1$-forms $($the cotangent sheaf$)$ $\Omega^1_X$ on a Poisson mani\-fold $X$ with Poisson bivector field $\pi \in H^0\bigl(X, \wedge^2 \mathcal{T}_X\bigr)$, has a canonical Lie algebroid structure, described as follows.
	The anchor map $\tilde{\pi}\colon \Omega^1_X \rightarrow \mathcal{T}_X$ is defined by $\tilde{\pi}(\omega) (f)= \pi(\omega, {\rm d}f)$ for $f \in \mathcal{O}_X$, $\omega \in \Omega^1_X$, and the Lie bracket on $\Omega^1_X$ is defined using Lie derivatives by
$[\omega, \omega']:= L_{\tilde{\pi}(\omega)} \omega' - L_{\tilde{\pi}(\omega')} \omega - {\rm d}(\pi(\omega, \omega'))$ for $\omega, \omega' \in \Omega^1_X$ (see \cite{RF, BP}).
\end{Exm}	
\begin{Exm}
	\label{analytic spaces}
	Let $X$ be a complex manifold and $\mathcal{O}_X$ denote the sheaf of holomorphic functions. The vanishing set $($or the zero locus$)$ $Y:= V(\mathcal{I})$ of an ideal-sheaf $\mathcal{I} \subset \mathcal{O}_X$ is a subspace of $X$, which is not necessarily a submanifold and may have singularities. The sheaf of functions $\mathcal{O}_Y := \mathcal{O}_X/{\mathcal{I}}$ on $Y$ serves as its structure sheaf, and the pair $(Y, \mathcal{O}_Y)$ constitutes an analytic space (see in \cite{BP}).
	
Here, the tangent sheaf $\mathcal{T}_Y := (\mathcal{D}{\rm er}_{\mathbb{C}_Y}(\mathcal{O}_Y), [\cdot,\cdot]_c, {\rm id}_Y)$ of $Y$ is not necessarily a locally free Lie algebroid over $(Y, \mathcal{O}_Y)$ (for example, if $Y$ is a normal crossing divisor). Consider the sheaf of logarithmic derivations $($forms a generalized involutive distribution or singular foliation$)$ \cite{AA, BP} as
\[
\mathcal{T}_X(-\log Y) := \{D \in \mathcal{T}_X \mid D(\mathcal{I}) \subset \mathcal{I}\} \hookrightarrow \mathcal{T}_X
\]
$($geometrically, it represents the sheaf of the vector fields on $X$ whose flow preserves the divisor~$Y)$ with the canonical Lie algebroid structure. It is associated with $\mathcal{T}_Y$ via the canonical map
\[
 \rho \colon \ \mathcal{T}_X(-\log Y) \rightarrow \mathcal{T}_Y
\]
	defined by $\rho (g D) = \bar{g} \tilde{D}$, for $g \in \mathcal{O}_X$, $D \in\mathcal{T}_X(-\log Y)$ as $\tilde{D}\bigl(\bar{f}\bigr) = \overline{D(f)}$, where $\bar{h}=h + \mathcal{I} \in \mathcal{O}_Y$ for $h \in \mathcal{O}_X$.
\end{Exm}

There exist examples of Lie algebroids, such as free Lie algebroids \cite{MK} and sheaf of primitive elements~\cite{AA}, where the underlying module structure is not a coherent sheaf but rather a quasi-coherent sheaf.

\subsection{Gerstenhaber algebras over a topological ringed space} \label{sheaves of G-algs}

In the study of deformations of associative algebras, M. Gerstenhaber discovered an algebraic structure \cite{MGer, MG-SS}, which is now known as the Gerstenhaber algebra. Gerstenhaber algebras also appears in the study of the space of multisections of a Lie algebroid over a smooth mani\-fold~\mbox{\cite{CV, PX}}. These notions appear in the study of the formality theory of Lie algebroids~\cite{DC}, which subsequently leads to a generalization of the theory of deformation quantization (see~\cite{TS}). We consider sheaves of Gerstenhaber algebras as follows.

\begin{Def} \label{sheaf of G-algs and morphisms} A sheaf of Gerstenhaber algebras on $X$ or a Gerstenhaber algebra over $(X,\mathcal{O}_X)$ consists of a triple $\bigl(\mathcal{A}:=\bigoplus_{i \in \mathbb{Z}} \mathcal{A}^i, \wedge, [\cdot, \cdot]\bigr)$, where each $\mathcal{A}^i$ are $\mathcal{O}_X$-modules satisfying the following conditions:
	\begin{itemize}\itemsep=0pt
		\item $(\mathcal{A}, \wedge)$ is a graded commutative associative $\mathcal{O}_X$-algebra.
		\item $(\mathcal{A}[1], [\cdot,\cdot])$ is a graded $\mathbb{K}_X$-Lie algebra, where $\mathcal{A}[1]:=\bigoplus_{i \in\mathbb{Z}} \mathcal{A}^{(i)}$ with $\mathcal{A}^{(i)}=\mathcal{A}^{i+1}$.
		\item For each $a \in \mathcal{A}^{(i)}$, the section $[a,\cdot]$ is a derivation with respect to $\wedge$ with degree $i$, i.e.,
		the operations $\wedge$ and $[\cdot,\cdot]$ are compatible through the following Leibniz rule: for $a\in \mathcal{A}^i$, $b\in \mathcal{A}^j$ and $c\in \mathcal{A}$,
\[
[a, b \wedge c]=[a,b] \wedge c + (-1)^{(i-1)j}b \wedge [a,c].
\]
\end{itemize}
	
	A morphism $\psi\colon \bigl(\mathcal{A}:=\bigoplus_{i \in \mathbb{Z}} \mathcal{A}^i, \wedge, [\cdot, \cdot]\bigr) \rightarrow \bigl(\mathcal{A}':=\bigoplus_{i \in \mathbb{Z}} \mathcal{A'}^i, \wedge', [\cdot, \cdot]'\bigr)$ of Gerstenhaber algebras over a topological ringed space $(X, \mathcal{O}_X)$ is given as follows:
	\begin{itemize}\itemsep=0pt
		\item The map $\psi\colon (\mathcal{A}, \wedge) \rightarrow (\mathcal{A}', \wedge')$ is a graded $\mathcal{O}_X$-algebra homomorphism, and
		\item the map $\psi\colon (\mathcal{A}[1], [\cdot, \cdot]) \rightarrow (\mathcal{A}'[1], [\cdot, \cdot]')$ is a graded $\mathbb{K}_X$-algebra homomorphism.
	\end{itemize}
\end{Def}

\begin{Rem}
$($Sheaves of$)$ Gerstenhaber algebras are generalization $($or graded analogue$)$ of $($sheaves of$)$ Poisson algebras.
\end{Rem}

\begin{Rem} \label{G-alg to Lie alg} In most of the cases, we consider $\mathcal{A}^i=0$ for $i<0$ and $\mathcal{A}^0=\mathcal{O}_X$ for a Gerstenhaber algebra $\bigl(\mathcal{A}:=\bigoplus_{i \in \mathbb{Z}} \mathcal{A}^i, \wedge, [\cdot, \cdot]\bigr)$ over $(X, \mathcal{O}_X)$.
	Notice that for such a Gerstenhaber algebra $\bigl(\mathcal{A} := \bigoplus_{i \geq 0} \mathcal{A}^i, \wedge, [\cdot, \cdot]\bigr)$, the $\mathcal{O}_X$-module $\mathcal{A}^1$ canonically forms a Lie algebroid over $(X, \mathcal{O}_X)$, as described follows. The restriction of the G-bracket gives the $\mathbb{K}_X$-Lie algebra $\bigl(\mathcal{A}^1, [\cdot, \cdot]|_{\mathcal{A}^1}\bigr)$, and restricting the multiplication $\wedge\colon \mathcal{O}_X \otimes_{\mathbb{K}_X} \mathcal{A}^1 \rightarrow \mathcal{A}^1$ gives the $\mathcal{O}_X$-module structure. The anchor map is given by the assignment $D \mapsto [D, \cdot]$ for $D \in \mathcal{A}^1$.
	
	This construction globalizes the algebraic counterpart discussed in \cite{MG-SS, JH2}.
\end{Rem}
\begin{Exm} \label{multivector fields}
	For a topological ringed space $(X, \mathcal{O}_X)$, the sheaf of multiderivations over~$\mathcal{O}_X$~is
\[
\wedge^\bullet_{\mathcal{O}_X}\mathcal{T}_X := \left(\bigoplus^{\infty}_{i=0} \wedge_{\mathcal{O}_X}^i \mathcal{T}_X, \wedge \right),
\]
 forms a sheaf of Gerstenhaber algebras. It is the sheafification of the presheaf $U \mapsto \wedge^\bullet_{\mathcal{O}_X(U)} \mathcal{T}_X(U)$ of exterior algebras for $\mathcal{T}_X(U)$ over $\mathcal{O}_X(U)$, where $U$ varies in open sets of~$X$.
	In particular, when $X$ is a non-singular space, this is equivalent to the sheaf of multivector fields $\wedge^\bullet_{\mathcal{O}_X}\mathfrak{X}_X$ on~$X$. (See Example~\ref{Example8}.)
	
	The (sheaf of) Gerstenhaber algebra structure is given by the Schouten--Nijenhuis bracket
\[
[\cdot,\cdot]\colon \
\wedge^\bullet_{\mathcal{O}_X}\mathcal{T}_X[1]
\otimes_{\mathcal{O}_X}
\wedge^\bullet_{\mathcal{O}_X}\mathcal{T}_X[1]
\rightarrow
\wedge^\bullet_{\mathcal{O}_X}\mathcal{T}_X[1],
\]
	which is defined as (on the sheaf of multisections)
\begin{gather*}
	[D_1\! \wedge \cdots \wedge D_m, D_{m+1}\! \wedge \cdots \wedge D_n] :=
\!\sum_{i \leq m < j}\!\! (-1)^{i+j+m}[D_i, D_j]_c \wedge D_1\! \wedge \cdots \hat{D}_i \cdots \hat{D}_j \cdots \wedge D_n,
\end{gather*}
	where $D_i \in \mathcal{T}_X$ for $i = 1, \dots, n$, and $[\cdot, \cdot]_c$ denotes the commutator Lie bracket on $\mathcal{T}_X$.
\end{Exm}
\begin{Exm} \label{Lie alg to G-alg}
	For a Lie algebroid $\mathcal{L}$ over $(X,\mathcal{O}_X)$, the standard Gerstenhaber algebra associated to $\mathcal{L}$ or the standard Gerstenhaber algebra structure on $\wedge_{\mathcal{O}_X}^\bullet \mathcal{L}$ is given by the triple \[
\mathcal{G}_{\mathcal{L}}:=\left(\bigoplus^{\infty}_{i=0} \wedge_{\mathcal{O}_X}^i \mathcal{L}, \wedge, [\cdot,\cdot]^{\mathcal{L}}\right),
\]
 where the graded $\mathbb{K}_X$-Lie bracket $[\cdot,\cdot]^{\mathcal{L}}$ is the generalization of the Schouten--Nijenhuis bracket in the context of Lie algebroids, replacing the commutator bracket $[\cdot,\cdot]_c$ of $\mathcal{T}_X$ by the Lie bracket $[\cdot,\cdot]_{\mathcal{L}}$ of $\mathcal{L}$.
\end{Exm}

In the following we describe an algebro-geometric counterpart of the relationship between Lie algebroids and Gerstenhaber algebras (see \cite{JH1, PX} for the smooth and algebraic case).
\begin{Thm} \label{G-alg & Lie alg}
	For an $\mathcal{O}_X$-module $\mathcal{E}$, there is a $($sheaf of$)$ Gerstenhaber algebra structure on the exterior algebra $\wedge_{\mathcal{O}_X}^\bullet \mathcal{E}$ if and only if $\mathcal{E}$ has a Lie algebroid structure over $(X, \mathcal{O}_X)$.
\end{Thm}
\subsection{Chevalley--Eilenberg--de Rham complex} \label{C-E-d cohomology}
The Chevalley--Eilenberg--de Rham complex is commonly used to study the cohomology of Lie algebroids over smooth manifolds \cite{KM}. Additionally, its algebro-geometric counterparts are discussed in \cite{UB, BRT, BP}. Here, we present a sheaf-theoretic generalization of this cochain complex, using the Lie--Rinehart cochain complex \cite{GR}.

For a Lie algebroid $\mathcal{L}:=(\mathcal{L}, [\cdot, \cdot], \mathfrak{a})$ over $(X,\mathcal{O}_X)$, first recall some of the key notions in the following.
Denote $^p {\rm Hom}_{\mathcal{O}_X(U)}\bigl(\wedge^\bullet_{\mathcal{O}_X(U)}\mathcal{L}(U), \mathcal{O}_X(U)\bigr)$ (or in short $^p \bigl(\wedge^\bullet_{\mathcal{O}_X(U)}\mathcal{L}(U)\bigr)^*$ by the
module of $\mathcal{O}_X(U)$-linear maps from $\wedge^\bullet_{\mathcal{O}_X(U)}\mathcal{L}(U)$ to $\mathcal{O}_X(U)$ which induces $\mathcal{O}_X(V)$-module homomorphisms from $\wedge^\bullet_{\mathcal{O}_X(V)}\mathcal{L}(V)$ to $\mathcal{O}_X(V)$ for each open subset $V$ in $U$,
and are compatible with restriction maps of the presheaf $U \mapsto \wedge^\bullet_{\mathcal{O}_X(U)}\mathcal{L}(U)$ and of the sheaf $\mathcal{O}_X$.

In the underlying process, we consider the presheaf $^p \wedge^{\bullet}_{\mathcal{O}_X} \mathcal{L}$ given by $U \mapsto \wedge^{\bullet}_{\mathcal{O}_X(U)} \mathcal{L}(U)$ together with all the presheaf homomorphisms from $^p \wedge^{\bullet}_{\mathcal{O}_X} \mathcal{L}$ to $\mathcal{O}_X$.

\subsubsection[The complex of a Lie algebroid L]{The complex of a Lie algebroid $\boldsymbol{\mathcal{L}}$}

Consider the cochain complex of $\mathcal{O}_X$-modules associated to a Lie algebroid $\mathcal{L}$
\begin{align} \label{Chevalley-Eilenberg-de Rham complex}
	\Omega^{\bullet}_{\mathcal{L}}:= \bigl(\mathcal{H}{\rm om}_{\mathcal{O}_X}\bigl(\wedge^\bullet_{\mathcal{O}_X}\mathcal{L}, \mathcal{O}_X\bigr), {\rm d}_{\mathcal{L}}\bigr).
\end{align}
It is the \emph{sheafification of the presheaf of cochain complexes} (see \cite{CW1} for such notions)
\[
U\mapsto \bigl({}^p {\rm Hom}_{\mathcal{O}_X(U)}\bigl(\wedge^\bullet_{\mathcal{O}_X(U)}\mathcal{L}(U), \mathcal{O}_X(U)\bigr), {\rm d}_{\mathcal{L}(U)}\bigr)=:\Omega_{\mathcal{L}(U)}^\bullet,
\]
and the differential ${\rm d}_{\mathcal{L}(U)}\colon {}^p \bigl(\wedge^\bullet_{\mathcal{O}_X(U)}\mathcal{L}(U)\bigr)^* \rightarrow {}^p \bigl(\wedge^{\bullet +1}_{\mathcal{O}_X(U)}\mathcal{L}(U)\bigr)^*$ is given by
\begin{align}
&{\rm d}_{\mathcal{L}(U)}(\omega)(D_1\wedge \cdots \wedge D_{k+1}) \nonumber\\
& \qquad = \sum^{k+1}_{i=1}(-1)^{i+1}\mathfrak{a}(U)(D_i)(\omega(D_1\wedge \cdots \wedge \hat{D_i} \wedge \cdots \wedge D_{k+1}))\nonumber\\
&\qquad\quad + \sum_{i< j}(-1)^{i+j}\omega([D_i, D_j]\wedge D_1\wedge \cdots \wedge \hat{D_i}\wedge \cdots \wedge \hat{D_j}\wedge \cdots \wedge D_{k+1}),\label{C-E-D differential}
\end{align}
where $D_1, \dots, D_{k+1} \in \mathcal{L}(U)$ and $ \omega \in \Omega^k_{\mathcal{L}(U)}$, associated with an open subset $U$ of $X$.

Equivalently, one can consider $\Omega^{\bullet}_{\mathcal{L}}(U)= \bigl({\rm Hom}_{\mathcal{O}_U}\bigl(\wedge^\bullet_{\mathcal{O}_U}\mathcal{L}_U, \mathcal{O}_U\bigr), {\rm d}_{\mathcal{L}}(U)\bigr)$ for an open set ${U \subset X}$, where a section in $\Omega^{\bullet}_{\mathcal{L}}(U)$ is a sheaf homomorphism ($\mathcal{O}_U$-linear map) from $\wedge^\bullet_{\mathcal{O}_U}\mathcal{L}_U$ to $\mathcal{O}_U$ compatible with restriction of the sheaf $V \mapsto {\rm Hom}_{\mathcal{O}_V}\bigl(\wedge^\bullet_{\mathcal{O}_V}\mathcal{L}_V, \mathcal{O}_V\bigr)$ and the differential $V \mapsto {\rm d}_{\mathcal{L}}(V)$, for open subset $V \subset U$.

Notice that, the differential ${\rm d}_{\mathcal{L}}$ is a $\mathbb{K}_X$-linear map but not an $\mathcal{O}_X$-module homomorphism, satisfies the graded Leibniz rule: for $\omega_1 \in \Omega^k_{\mathcal{L}}$ and $\omega_2 \in \Omega^l_{\mathcal{L}}$, $k, l \in \mathbb{N} \cup \{0\}$, we obtain
\[
{\rm d}_\mathcal{L}(\omega_1 \wedge \omega_2)= {\rm d}_\mathcal{L} \omega_1 \wedge \omega_2 + (-1)^k \omega_1 \wedge {\rm d}_\mathcal{L}\omega_2.
\]

Note that if $\mathcal{L}$ is a locally free $\mathcal{O}_X$-modules of finite rank, then we have
\[
\mathcal{H}{\rm om}_{\mathcal{O}_X}\bigl(\wedge^{\bullet}_{\mathcal{O}_X}\mathcal{L}, \mathcal{O}_X\bigr) \cong \wedge^{\bullet}_{\mathcal{O}_X}\mathcal{L}^*.
\]
This complex is equivalent to the standard complex associated with a smooth Lie algebroid~\cite{KM} or a holomorphic Lie algebroid~\cite{BP}. In particular, we get the Chevalley-Eilenberg complex when~$\mathcal{L}$ is a Lie algebroid over a point; and the de Rham complex $\Omega^\bullet_X$ when $\mathcal{L}=\mathcal{T}_X$ over a~nonsingular space $X$. In general, the complex $\Omega^\bullet_X$ over $(X, \mathcal{O}_X)$ produces \emph{Zariski differentials}, not necessarily \emph{K\"ahler differentials} (see \cite{JH, CK}).

\begin{Rem}
Moreover, one can consider the cochain complex of K\"ahler differentials over a~scheme $(X, \mathcal{O}_X)$, given by
$\bigl(\wedge^{\bullet}_{\mathcal{O}_X} \Omega^1_X, {\rm d}_X\bigr)$,
	where $\Omega^1_X$ is the $\mathcal{O}_X$-module of K\"ahler differential 1\nobreakdash-forms on $X$, and the differential ${\rm d}_X$ extends the canonical map ${\rm d}^1_X\colon \mathcal{O}_X \to \Omega^1_X$.
	If $\Omega^1_X$ is reflexive, then $(\mathcal{T}_X)^* = \bigl(\Omega^1_X\bigr)^{**} \cong \Omega^1_X$ as $\mathcal{O}_X$-modules, a property that holds when $X$ is a~smooth scheme of finite type. Consequently, the algebraic de Rham complex of $X$ is recovered by substituting $\mathcal{L} = \mathcal{T}_X$ in the Lie algebroid complex discussed in~\cite{UB}.
\end{Rem}
Next, we consider an algebro-geometric analogue of the correspondence between Lie algebroid structure on a smooth vector bundle and differential graded algebra $($DGA$)$ structure on the exterior algebra of its dual bundle (see \cite{KS-M, PX}).
\begin{Thm} \label{Lie alg & dga} Given a locally free $\mathcal{O}_X$-module $\mathcal{E}$ of finite rank, $\mathcal{E}$ has a Lie algebroid structure over $(X, \mathcal{O}_X)$ if and only if the exterior algebra of sheaf of multisections \smash{$\wedge^{\bullet}_{\mathcal{O}_X} \mathcal{E}^*:=\bigoplus^{\infty}_{k=0}\wedge^k_{\mathcal{O}_X} \mathcal{E}^*$} of $\mathcal{E}^*=\mathcal{H}{\rm om}_{\mathcal{O}_X}(\mathcal{E}, \mathcal{O}_X)$ forms a sheaf of differential graded algebras over $(X, \mathcal{O}_X)$.
\end{Thm}

First, we recall some facts about \emph{hypercohomology} (see \cite{UB, VID, CW1}) in order to consider an analogue of the Chevalley--Eilenberg--de Rham cohomology.

\subsubsection[Hypercohomology of the complex Omega\^{}bullet\_L]{Hypercohomology of the complex $\boldsymbol{\Omega^{\bullet}_{\mathcal{L}}}$} \label{Hypercohomology}

We begin by recalling the hypercohomology for a cochain complex of sheaves of modules. We then describe the hypercohomology of the complex $\Omega^{\bullet}_{\mathcal{L}}$.

Let $(X, \mathcal{O}_X)$ be a topological ringed space with $\mathcal{O}_X$-modules $\mathcal{C}^i$, for $i \geq 0$. Consider a cochain complex $(\mathcal{C}^\bullet, {\rm d})$ where $\mathcal{C}^\bullet:= \bigoplus_i \mathcal{C}^i$ is a graded $\mathcal{O}_X$-module and ${\rm d}\colon \mathcal{C}^\bullet \rightarrow \mathcal{C}^{\bullet +1}$ is the co-boundary map, i.e., ${\rm d}$ is a $\mathbb{K}_X$-linear map of degree $1$ satisfying ${\rm d}^2=0$. Similarly, for a chain complex, instead of a co-boundary map we have a boundary map $\partial\colon \mathcal{C}^\bullet \rightarrow \mathcal{C}^{\bullet -1}$, a $\mathbb{K}_X$-linear map of degree $-1$ satisfying $\partial^2=0$.

To address the following topics of discussion consistently, we consider such (co)chain complexes where $\mathcal{O}_X$-linearity of the (co)boundary maps is not necessarily required.

For an $\mathcal{O}_X$-module $\mathcal{F}$ on $X$, there exists an injective resolution $\mathcal{F} \overset{\sim}{\rightarrow} \mathcal{C}^\bullet(\mathcal{F})$ $($a quasi-isomorphism) of $\mathcal{O}_X$-modules, known as the flabby Godement resolution of $\mathcal{F}$. Thus, for a~complex of~$\mathcal{O}_X$-modules $\mathcal{F}^\bullet$, we consider the associated bicomplex of sheaves $\mathcal{C}^\bullet(\mathcal{F}^\bullet)=(\mathcal{C}^p(\mathcal{F}^q))$ $(p,q \in \mathbb{N} \cup \{0\})$ of injective resolutions. The original complex is embedded in the total complex $\mathcal{K}^\bullet= {\rm tot}(\mathcal{C}^\bullet(\mathcal{F}^\bullet))$, and this embedding is a quasi-isomorphism. The cohomology of the associated complex of global sections $\mathcal{K}^\bullet(X)={\rm tot}(\mathcal{C}^\bullet(\mathcal{F}^\bullet))(X)$ is called the hypercohomology of $\mathcal{F}^\bullet$, and denoted by $\mathbb{H}^\bullet(X, \mathcal{F}^\bullet)$.

Denote the category of cochain complexes of sheaves of $\mathcal{O}_X$-modules over $(X,\mathcal{O}_X)$ by $\operatorname{Shv}(X)^\bullet$ and the category of $\mathbb{K}$-vector spaces by $\operatorname{Vect}_{\mathbb{K}}$.
For $\mathcal{F}^\bullet\in \operatorname{Shv}(X)^\bullet$, the $k$-th cohomology sheaf of the complex $\mathcal{F}^\bullet$ is
\[
\mathcal{H}^k(\mathcal{F}^\bullet) := \mathcal{K}{\rm er} \bigl(\mathcal{F}^k \rightarrow \mathcal{F}^{k+1}\bigr)/{\mathcal{I}{\rm m}\bigl(\mathcal{F}^{k-1} \rightarrow \mathcal{F}^{k}\bigr)}
\]
$($where $\mathcal{F}^{-1}:=\{0\})$ for $k \in \mathbb{N} \cup \{0\}$, considered as quotient of $\mathbb{K}_X$-vector spaces.
Moreover, a map of complexes $\mathcal{F}^\bullet \rightarrow \mathcal{G}^\bullet$ is a quasi-isomorphism if the induced map on the cohomology sheaves $\mathcal{H}^k(\mathcal{F}^\bullet) \rightarrow \mathcal{H}^k(\mathcal{G}^\bullet)$ is an isomorphism for all~$k$.
The $k$-th hypercohomology is a functor
\[
\mathbb{H}^k(X,-)\colon \ \operatorname{Shv}(X)^\bullet \rightarrow \operatorname{Vect}_{\mathbb{K}}
\]
 for $k \in \mathbb{N} \cup \{0\}$, that satisfies the following two conditions:
\begin{itemize}\itemsep=0pt
	\item A quasi-isomorphism of complexes $f^\bullet\colon \mathcal{F}^\bullet \rightarrow \mathcal{G}^\bullet$ induces an isomorphism $\mathbb{H}^k(X, f^\bullet)$.
	\item
	If $\mathcal{I}^\bullet$ is a complex of injective sheaves, then $\mathbb{H}^k(X, \mathcal{I}^\bullet) = H^k(\Gamma(X, \mathcal{I}^\bullet))$ or $H^k(\mathcal{I}^\bullet (X))$.
\end{itemize}

Denote the graded vector space \smash{$\bigoplus_{n=0}^{\infty} \mathbb{H}^n(X,\mathcal{F}^\bullet)$} as $\mathbb{H}^\bullet(X,\mathcal{F}^\bullet)$.

The hypercohomology of the cochain complex (\ref{Chevalley-Eilenberg-de Rham complex}) of $\mathcal{O}_X$-modules is called the Lie algebroid cohomology or the \emph{Chevalley--Eilenberg--de Rham cohomology} of $\mathcal{L}$ with trivial coefficient and is denoted by $\mathbb{H}^\bullet(\mathcal{L}, \mathcal{O}_X)$.

\begin{Rem} \label{hyperhomology}
Later, we consider an analogue of the Lie algebroid homology \cite{JH1, PX}, which is obtained from a similar version of the above description of hypercohomology in the context of chain complexes of sheaves (see~\cite{CW1}), which we call the hyperhomology of the Lie algebroid.
\end{Rem}

\subsection{Lie bialgebroids over a topological ringed space} \label{Lie bialgebroids}
Lie bialgebroids over a smooth manifold, consists of two compatible Lie algebroids defined on dual vector bundles, is a standard notion appears in the study of Poisson geometry \cite{KS,KM,KM-PX,PX}.
Here, we consider an analogue of it in the algebro-geometric settings.
\begin{Def} A Lie bialgebroid over a topological ringed space $(X,\mathcal{O}_X)$ is a pair $(\mathcal{L}, \mathcal{L}^*)$ of locally free Lie algebroids of finite rank in duality, where the Lie brackets satisfy a condition which can be expressed in terms of the differential ${\rm d}_{\mathcal{L}^*}=:{\rm d}_*$ on $\wedge^\bullet_{\mathcal{O}_X} \mathcal{L}$ defined by the Lie algebroid structure of $\mathcal{L}^*$ and the Gerstenhaber bracket $[\cdot, \cdot]^{\mathcal{L}}$ on $\wedge^\bullet_{\mathcal{O}_X}\mathcal{L}$ defined by the Lie algebroid structure on $\mathcal{L}$, is given by
\[
{\rm d}_*[D_1, D_2]= [{\rm d}_*D_1, D_2]^{\mathcal{L}} + [D_1, {\rm d}_*D_2]^{\mathcal{L}}
\]
for all sections $D_1$, $D_2 $ of $\mathcal{L}$, where $[\cdot, \cdot]:=[\cdot, \cdot]_{\mathcal{L}}=[\cdot, \cdot]^{\mathcal{L}}|_{\mathcal{L}}$ is the Lie bracket on~$\mathcal{L}$.
	
In other words,	a Lie bialgebroid over a topological ringed space $(X,\mathcal{O}_X)$ is a dual pair $(\mathcal{L}, \mathcal{L}^*)$ of locally free Lie algebroids of finite rank such that the differential~${\rm d}_*$, induced from the Lie algebroid structure on $\mathcal{L}^*$, is a derivation of the graded Lie bracket $[\cdot, \cdot]^{\mathcal{L}}$ on $\wedge^\bullet_{\mathcal{O}_X}\mathcal{L}$.
\end{Def}
\begin{Rem}
	The condition $\mathcal{L}^*$ is a locally free $\mathcal{O}_X$-modules of finite rank implies
\[
\mathcal{H}{\rm om}_{\mathcal{O}_X}\bigl(\wedge^{\bullet}_{\mathcal{O}_X}\mathcal{L}^*, \mathcal{O}_X\bigr) \cong \wedge^{\bullet}_{\mathcal{O}_X}\mathcal{L}.
\]
\end{Rem}

\begin{Exm}
	For a Poisson manifold $(X, \pi)$, the dual pair $\bigl(\mathcal{T}_X, \Omega^1_X\bigr)$ forms a Lie bialgebroid, where the Lie bracket on the cotangent sheaf $\Omega^1_X$ and the anchor map $\tilde{\pi}\colon \Omega^1_X \rightarrow \mathcal{T}_X$ is given by the Poisson bivector field $\pi$ on $X$ $($or, by the induced Poisson bracket $\{\cdot, \cdot \}$ on $\mathcal{O}_X)$, see Example~\ref{Cotantgent Poisson}. More generally, there is an analogous notion, known as triangular Lie bialgebroid in the context of Lie algebroids \cite{KS}.
\end{Exm}
\begin{Rem}
	Let $(\mathcal{L}, \mathcal{L}^*)$ be a pair of locally free Lie algebroids of finite rank in duality over a topological ringed space. Then
	$(\mathcal{L}, \mathcal{L}^*)$ is a Lie bialgebroid if and only if $(\mathcal{L}^*, \mathcal{L})$ is a Lie bialgebroid.
\end{Rem}
\subsubsection{(Sheaves of) Strong differential Gerstenhaber algebras} \label{strong diff G-algs}
Lie bialgebroids over smooth manifolds and strong differential Gerstenhaber algebras are closely related \cite{KS, PX}. Here, we define strong differential Gerstenhaber algebras over topological ringed spaces and extend this relationship.
\begin{Def}
	A \emph{strong differential Gerstenhaber algebra over a topological ringed space} $(X, \mathcal{O}_X)$ is a Gerstenhaber algebra $\mathcal{A}^\bullet:=\bigl(\bigoplus^{\infty}_{i=0} \mathcal{A}^i, \wedge, [\cdot, \cdot]\bigr)$ over $(X,\mathcal{O}_X)$ equipped with a~differential $($a~$\mathbb{K}_X$-linear map of degree $1)$
	$\tilde{{\rm d}}\colon \bigoplus^{\infty}_{i=0} \mathcal{A}^i \rightarrow \bigoplus^{\infty}_{i=0} \mathcal{A}^{i+1}$, which is a derivation with respect to $\wedge$ satisfying $\tilde{{\rm d}}^2=0$, and it is also a derivation of the graded Lie bracket $[\cdot, \cdot]$.
\end{Def}

An algebro-geometric analogue of the correspondence between these notions~\cite{KS} is described as follows.
\begin{Thm} \label{SDGA & LBA}
For a locally free $\mathcal{O}_X$-module $\mathcal{E}$ of finite rank over a topological ringed space $(X,\mathcal{O}_X)$, the exterior algebra $\wedge_{\mathcal{O}_X}^\bullet\mathcal{E}$ has a strong differential Gerstenhaber algebra structure over $(X,\mathcal{O}_X)$ if and only if the dual pair $(\mathcal{E}, \mathcal{E}^*)$ has a Lie bialgebroid structure over $(X,\mathcal{O}_X)$.
\end{Thm}

\subsection{(Sheaves of) BV-algebras or exact Gerstenhaber algebras}
BV-algebras play an important role in geometry and mathematical physics (see \cite{JH1, KS, PX}). We consider an algebro-geometric analogue of this structure.
\begin{Def}	An exact Gerstenhaber algebra (or BV-algebra) over a topological ringed space $(X,\mathcal{O}_X)$ is a Gerstenhaber algebra $\bigl(\mathcal{A}:=\bigoplus^{\infty}_{i=0} \mathcal{A}^i, \wedge, [\cdot, \cdot]\bigr)$ over $(X,\mathcal{O}_X)$ equipped with an operator $\partial$ of degree $-1$, square zero, and generating the Gerstenhaber bracket as
\[
[a, b]= (-1)^{i}(\partial(a \wedge b)- \partial a \wedge b- (-1)^i a \wedge \partial b ),
\]
	where $a \in \mathcal{A}^i \subset \mathcal{A}$ and $b \in \mathcal{A}$. The operator $\partial$ is called a generating operator.	
\end{Def}
We recall flat connections of a Lie algebroid using \emph{Atiyah algebroids} \cite{UB,MK,BP,AS,PT}.

\subsubsection{Atiyah algebroids and Lie algebroid connections} \label{Atiyah and Lie connections}
For an $\mathcal{O}_X$-module $\mathcal{E}$, we form a Lie algebroid consisting of the sheaf of differential operators on~$\mathcal{E}$ of order $\leq 1$ with scalar symbols, i.e.,
\begin{gather*}
\mathcal{A}{\rm t}(\mathcal{E}):= \{D\in \mathcal{E}{\rm nd}_{\mathbb{K}_X}(\mathcal{E})\mid
D(fs)=fD(s)+\sigma_D(f)s \ \text{for a unique} \ \sigma_D\in \mathcal{T}_X, \\
\hphantom{\mathcal{A}{\rm t}(\mathcal{E}):= \{D\in \mathcal{E}{\rm nd}_{\mathbb{K}_X}(\mathcal{E})\mid{}}{} \text{where} \ f \in \mathcal{O}_X,\, s\in \mathcal{E}\},
\end{gather*}
(thus, here $\sigma_D(f)=[D, f]_c \in \mathcal{O}_X$ for $D \in \mathcal{A}{\rm t}(\mathcal{E})$ and $f \in \mathcal{O}_X$ holds) with the anchor map defined by
\[
	 \sigma\colon \ \mathcal{A}{\rm t}(\mathcal{E}) \rightarrow \mathcal{T}_X,\qquad \text{where $D \mapsto \sigma_D$}
\]
and the Lie bracket is commutator bracket. This Lie algebroid is the so-called Atiyah algebroid of the $\mathcal{O}_X$-module $\mathcal{E}$. It provides a short exact sequence (s.e.s.) of Lie algebroids over $(X, \mathcal{O}_X)$
(an abelian Lie algebroid extension),
\begin{align}\label{Lie algebroid extension}
	0 \rightarrow \mathcal{E}{\rm nd}_{\mathcal{O}_X}(\mathcal{E}) \hookrightarrow \mathcal{A}{\rm t}(\mathcal{E}) \overset{\sigma}{\rightarrow} \mathcal{T}_X \rightarrow 0.
\end{align}

In particular, when $\mathcal{E}=\mathcal{O}_X$ for a nonsingular space $X$, we have $\mathcal{A}{\rm t}(\mathcal{O}_X)\cong \mathcal{O}_X \oplus \mathcal{T}_X$.

\begin{Def} \label{L-connections}
Given Lie algebroid $(\mathcal{L}, [\cdot, \cdot], \mathfrak{a})$ on a topological ringed space $(X, \mathcal{O}_X )$, an $\mathcal{L}$-connection on an $\mathcal{O}_X$-module $\mathcal{E}$ is defined by an $\mathcal{O}_X$-linear map
	\begin{align}\label{L-connection}
		\nabla\colon \ \mathcal{L}\rightarrow \mathcal{A}{\rm t}(\mathcal{E}) \qquad \text{given by} \ D\mapsto \nabla_D
	\end{align}
	satisfying the Leibniz rule $\nabla_D(f s)= f \nabla_D(s)+\mathfrak{a}(D)(f) s$ for sections $f \in \mathcal{O}_X$, $D \in \mathcal{L}$ and $s \in \mathcal{E}$ (see \cite{UB, BRT, MK, BP}). To ensure compatibility, we choose the condition $\sigma \circ \nabla= \mathfrak{a}$.
\end{Def}
\begin{Rem} \label{L-connection using dual}
	Equivalently, an $\mathcal{L}$-connection on $\mathcal{E}$ is described by a $\mathbb{K}_X$-linear map, if in addition $\mathcal{L}$ is a locally free $\mathcal{O}_X$-module of finite rank, as follows (see \cite{BRT, BP}):
\[
{\rm d}_{\nabla} \colon \ \mathcal{E} \rightarrow \Omega_{\mathcal{L}}^1 \otimes_{\mathcal{O}_X}\mathcal{E}
\] satisfying the Leibniz rule ${\rm d}_{\nabla}(fs)=f {\rm d}_{\nabla} s + {\mathfrak{a}^*}({\rm d}f)\otimes s$, where
	${\mathfrak{a}^*}\colon \Omega_X^1 \rightarrow \Omega_{\mathcal{L}}^1$ is the dual of the anchor map.
\end{Rem}

\begin{Def} \label{flat}
	An $\mathcal{L}$-connection on $\mathcal{E}$ is said to be flat if the map (\ref{L-connection}) is a Lie algebroid homomorphism, i.e., the map (\ref{L-connection}) satisfies curvature zero condition: for $D,D' \in \mathcal{L}$
	\begin{align*}
		R_{\nabla}(D,D'):=\nabla_{[D,D']}-[\nabla_D, \nabla_{D'}]_c=0.
	\end{align*}
	In that case, $(\mathcal{E},\nabla)$ is said to be a representation of $\mathcal{L}$ or an $\mathcal{L}$-module.
\end{Def}

\subsubsection[Particular cases of L-connections]{Particular cases of $\boldsymbol{\mathcal{L}}$-connections}

A connection $\nabla\colon \mathcal{T}_X \rightarrow \mathcal{A}{\rm t}(\mathcal{E})$ that satisfies $\sigma \circ \nabla=\mathrm{Id}_{\mathcal{T}_X}$ provides a splitting of the s.e.s.\ (\ref{Lie algebroid extension}) as $\mathcal{O}_X$-modules. If its curvature is zero then the s.e.s. splits as Lie algebroids.

Next, we consider some special cases where $X$ is a nonsingular space:
\begin{enumerate}\itemsep=0pt
\item[$(i)$] For $\mathcal{L}=\mathcal{T}_X$ and $\mathcal{E}$ a locally free $\mathcal{O}_X$-module of finite rank, a $\mathcal{T}_X$-connection on $\mathcal{E}$ always exists if $X$ is either a smooth manifold, or a Stein manifold, or an affine variety (see \cite{MA}).
\item[$(ii)$] For $\mathcal{L}=\mathcal{T}_X(- \log Y)$, if a connection exists on an $\mathcal{O}_X$-module $\mathcal{E}$, it is called a logarithmic connection on $\mathcal{E}$. It forms a meromorphic connection with simple poles along the divisor~$Y$ (see~\cite{CMD,BP}). For instance, if $X$ is a non-singular projective toric variety and $Y$ is the boundary divisor, then a locally free $\mathcal{O}_X$-module $\mathcal{E}$ admits a logarithmic connection if and only if it admits a torus equivariant structure (see \cite{IBMP, IB}).
\end{enumerate}

An algebro-geometric analogue of the correspondence between these notions \cite{PX} is described as follows.
\begin{Thm} \label{BV-alg & flat conn}
	Let $\mathcal{L}$ be a locally free $\mathcal{O}_X$-module of rank $n$. Then there is an exact Gerstenhaber algebra $($or BV-algebra$)$ structure on the exterior algebra $\wedge^{\bullet}_{\mathcal{O}_X}\mathcal{L}$ if and only if there is a~Lie algebroid structure on $\mathcal{L}$ with a flat $\mathcal{L}$-connection on the $\mathcal{O}_X$-module $\wedge^{n}_{\mathcal{O}_X}\mathcal{L}$.
\end{Thm}
\begin{Rem}\label{canonical BV-alg}
	Let $X$ be either a smooth manifold, or a Stein manifold, or a smooth affine variety, and $\mathcal{O}_X$ be the associated structure sheaf (see \cite{MA, FF}).
	For a locally free Lie algebroid $\mathcal{L}$ of rank $n$ for some $n \in \mathbb{N}$, there always exists a flat $\mathcal{L}$-connection on $\wedge^{n}_{\mathcal{O}_X}\mathcal{L}$. Thus, there always exists a canonical BV-algebra structure on $\wedge^{\bullet}_{\mathcal{O}_X}\mathcal{L}$ (see \cite{PX} for the smooth Lie algebroids case).
\end{Rem}

\section[L-Lie algebroids and G\_L-Gerstenhaber algebras]{$\boldsymbol{\mathcal{L}}$-Lie algebroids and $\boldsymbol{\mathcal{G}_{\mathcal{L}}}$-Gerstenhaber algebras}
\label{Sec 3}

We define $\mathcal{L}$-Lie algebroids for a given Lie algebroid $\mathcal{L}$ over a topological ringed space $(X, \mathcal{O}_X)$. To generalize the relationship between Lie algebroids and sheaves of Gerstenhaber algebras $($Theorem~\ref{G-alg & Lie alg}$)$ in this context, we introduce $\mathcal{A}$-Gerstenhaber algebras associated with a given Gerstenhaber algebra $\mathcal{A}$ over $(X, \mathcal{O}_X)$.

\begin{Def} \label{L-Lie algebroids} Let $(\mathcal{L}, [\cdot,\cdot]_{\mathcal{L}}, \mathfrak{a})$ be a Lie algebroid over $(X, \mathcal{O}_X)$. A Lie algebroid $(\mathcal{L'}, [\cdot,\cdot]_{\mathcal{L'}}, \mathfrak{a'})$ over $(X, \mathcal{O}_X)$ is said to be an $\mathcal{L}$-Lie algebroid if there exists a Lie algebroid homomorphism
\[
\phi\colon \ (\mathcal{L'}, [\cdot,\cdot]_{\mathcal{L'}}, \mathfrak{a'}) \rightarrow (\mathcal{L}, [\cdot,\cdot]_{\mathcal{L}}, \mathfrak{a}).
\]
	
	We denote this $\mathcal{L}$-Lie algebroid as $(\mathcal{L'},\phi)$.
\end{Def}

\begin{Exm}\quad
\begin{enumerate}\itemsep=0pt
	\item Consider the tangent sheaf $\mathcal{T}_X$ of a special ringed space $(X,\mathcal{O}_X)$ (see Remark~\ref{special spaces}). Any Lie algebroid $\mathcal{L}$ over $(X,\mathcal{O}_X)$ is a $\mathcal{T}_X$-Lie algebroid. For example, the cotangent sheaf over a Poisson manifold.
	\item Any Lie subalgebroid $\mathcal{L}'$ of a Lie algebroid $\mathcal{L}$ naturally inherits the structure of an $\mathcal{L}$-Lie algebroid via the inclusion map $\iota\colon \mathcal{L}' \to \mathcal{L}$. Singular foliations offer such examples, including the sheaf of logarithmic derivations and the characteristic (or orbit) foliation (see \cite{RF, BP}).
	
	\item For a triangular Lie bialgebroid $(\mathcal{L}^*, \mathcal{L}, \mathcal{P})$ where $\mathcal{P}$ is the Poisson $\mathcal{L}$-bivector, can be viewed as a Lie algebroid homomorphism $\mathcal{P}\colon \mathcal{L}^* \rightarrow \mathcal{L}$ (see \cite{KS}). Therefore, $(\mathcal{L}^*, \mathcal{P})$ is an $\mathcal{L}$-Lie algebroid.
\end{enumerate}
\end{Exm}

\begin{Exm}
Analogous to singular foliations (see Example~\ref{Foliation}), we have singular subalgebroids of a smooth Lie algebroid~$L$~\cite{MZ}.
An $L$-valued almost Lie algebroid $(L', \phi)$ over a singular subalgebroid $\mathcal{F}$ is
an $L$-Lie algebroid such that $\mathcal{I}{\rm m}(\Gamma(\phi))=\phi(\Gamma(L')) \subset \mathcal{F}$ of the Lie algebroid morphism $\phi\colon L' \rightarrow L$, where $\Gamma$ represents the space of global sections functor~\cite{RL25}. Considering the associated sheaf of sections, one gets an example of $\mathcal{L}$-Lie algebroids, where $\mathcal{L}=\Gamma_X(L)$.
\end{Exm}

\begin{Exm} \label{Lie groups}
Let $G$ be a (smooth, holomorphic, or algebraic) Lie group acting on a $G$-mani\-fold~$X$, and let $\bar{\phi}\colon G' \to G$ be a Lie group homomorphism for another Lie group $G'$. This map $\bar{\phi}$ induces a natural $G'$-action on the $G$-manifold $X$.
	Denote the Lie algebra of a Lie group $H$ by $\operatorname{Lie}(H)$, and the associated constant sheaf by $\operatorname{Lie}(H)_X$. Consider the associated transformation Lie algebroids $X \times \operatorname{Lie}(G')$ and $X \times \operatorname{Lie}(G)$ over $X$ (see~\cite{RF}) along with the canonical connecting map \smash{$\tilde{\phi}\colon X \times \operatorname{Lie}(G') \to X \times \operatorname{Lie}(G)$}.
	Using the $G$-action on~$X$, the map \smash{$\tilde{\phi}$} induces a Lie algebroid homomorphism over $(X, \mathcal{O}_X)$ at the level of sheaves of sections. Consequently, this defines an $\mathcal{O}_X \otimes_{\mathbb{K}_X} \operatorname{Lie}(G)_X$-Lie algebroid structure on~$\bigl(\mathcal{O}_X \otimes_{\mathbb{K}_X} \operatorname{Lie}(G')_X, {\rm id} \otimes \operatorname{Lie}\bigl(\bar{\phi}\bigr)\bigr)$.
\end{Exm}
\begin{Exm} \label{log Atiyah}
	The log Atiyah algebroid has appeared in the study of logarithmic connections on equivariant principal bundles (see \cite{IBMP, DKP}) or on equivariant vector bundles \cite{IB}. Here, we recall the notion log Atiyah algebroid $\mathcal{A}{\rm t}(\mathcal{E})(-\log Y)$ of an $\mathcal{O}_X$-module $\mathcal{E}$ associated with an analytic subspace $Y$ of a complex manifold $X$ (or, a divisor $Y$ of a smooth algebraic variety~$X)$,
	and describe it as a $\mathcal{T}_X(-\log Y)$-Lie algebroid. First, consider the Atiyah algebroid $\mathcal{A}{\rm t}({\mathcal{E}})$ of $\mathcal{E}$ and the Lie algebroid homomorphism $\iota\colon \mathcal{T}_X(-\log Y) \hookrightarrow \mathcal{T}_X$. The pull back of the Lie algebroid $\mathcal{T}_X(-\log Y)$ by the anchor map $\sigma\colon \mathcal{A}{\rm t}({\mathcal{E}}) \rightarrow \mathcal{T}_X$ is the log Atiyah algebroid of $\mathcal{E}$, i.e., $\mathcal{A}{\rm t}(\mathcal{E})(-\log Y) = \sigma^{-1}(\mathcal{T}_X(-\log Y))$. Thus, $(\mathcal{A}{\rm t}(\mathcal{E})(-\log Y), \bar{\sigma})$ is a $\mathcal{T}_X(-\log Y)$-Lie algebroid, where~$\bar{\sigma}$ is the restriction of the Lie algebroid homomorphism~$\sigma$ on $\mathcal{A}{\rm t}(\mathcal{E})(-\log Y)$.
\end{Exm}

\begin{Exm} \label{extensions of L}
	Given a Lie algebroid $(\mathcal{L}, [\cdot, \cdot], \mathfrak{a})$ over $(X, \mathcal{O}_X)$ and an element $[\omega] \in \mathbb{H}^2(\mathcal{L}, \mathcal{O}_X)$, we get a Lie algebroid extension of $\mathcal{L}$ by $\mathcal{O}_X$ (up to isomorphism) as follows (see \cite{HOM, PT2}).
	
Consider the $\mathcal{O}_X$-module $\mathcal{O}_X \oplus \mathcal{L}$ with the $\mathbb{K}_X$-Lie algebra structure given by
\[
[(f,D), (f',D')]_{\omega}:=(\mathfrak{a}(D)(f')-\mathfrak{a}(D')(f)+\omega(D,D'), [D, D'])
\]
	for $f, f'\in \mathcal{O}_X$, $D, D' \in \mathcal{L}$, together with the map
$
\mathfrak{a}_{\omega}\colon \mathcal{O}_X \oplus \mathcal{L} \rightarrow \mathcal{T}_X$, $(f, D) \mapsto \mathfrak{a}(D)$.
	Thus, $\mathcal{L}_{\omega}:=(\mathcal{O}_X \oplus \mathcal{L}, [\cdot, \cdot]_{\omega}, \mathfrak{a}_{\omega})$ forms a~Lie algebroid over $(X, \mathcal{O}_X)$, provides the abelian Lie algebroid extension (an equivalence class)
	$\mathcal{O}_X \hookrightarrow \mathcal{L}_{\omega} \twoheadrightarrow \mathcal{L}$.
	Thus, $(\mathcal{L}_{\omega}, \operatorname{pr}_2)$ is an $\mathcal{L}$-Lie algebroid for any $2$-cocycle $\omega$ of the complex $\Omega^{\bullet}_{\mathcal{L}}$ (see~(\ref{Chevalley-Eilenberg-de Rham complex})),
	where $\operatorname{pr}_2\colon \mathcal{O}_X \oplus \mathcal{L} \rightarrow \mathcal{L}$ denotes the projection onto the second component.
\end{Exm}

 \begin{Exm} \label{Free Lie algebroid} In~\cite{MK}, M.~Kapranov introduced the free Lie algebroid $\mathcal{FL}(\mathcal{E})$ of an anchored $\mathcal{O}_X$-module $\mathcal{E}$, in order to study the path algebroid $\mathcal{P}_X:=\mathcal{FL}(\mathcal{T}_X)$ of $X$. An anchored $\mathcal{O}_X$-module $(\mathcal{E}, \alpha)$ is an $\mathcal{O}_X$ module $\mathcal{E}$ together with an $\mathcal{O}_X$ module homomorphism $\alpha\colon \mathcal{E} \rightarrow \mathcal{T}_X$. A pre-module over $\mathcal{E}$ is an $\mathcal{O}_X$-module $\mathcal{F}$ together with an anchored module homomorphism $\mathcal{E} \rightarrow \mathcal{A}{\rm t}(\mathcal{F}).$
An important consequence is that pre-modules over $\mathcal{E}$ naturally form modules over the free Lie algebroid $\mathcal{FL}(\mathcal{E})$, follows from the universal property of $\mathcal{FL}(\mathcal{E})$. In particular, for $\mathcal{E}= \mathcal{T}_X$, a pre-module $\mathcal{F}$ is a connection (not necessarily flat) on $\mathcal{F}$.
	
Let $\mathcal{L}$ be a Lie algebroid whose anchor is $\mathfrak{a}$. For an anchored $\mathcal{O}_X$-module $(\mathcal{E}, \alpha)$ together with anchored module homomorphism $\phi\colon \mathcal{E} \rightarrow \mathcal{L}$, i.e., $\mathfrak{a} \circ \phi=\alpha$, we get the following. From \cite[Theorem~2.1.2]{MK},
	$\mathcal{FL}(\mathcal{E})$ has an $\mathcal{L}$-Lie algebroid structure if and only if there is an anchored module homomorphism from $\mathcal{E}$ to~$\mathcal{L}$,
\[
\xymatrix@C=4em@R=3em{
	{\mathcal{E}} \ar[r]^{\phi} \ar[d]_{\iota} &
	{\mathcal{L}} \ar[d]^{\mathfrak{a}} \\
	{\mathcal{FL}(\mathcal{E})} \ar[r]_{\psi} \ar@{.>}[ur] &
	{\mathcal{T}_X}.
}
\]
\end{Exm}
\begin{Rem} \label{Lie-Atiyah alg}
	Consider the Atiyah algebroid $\mathcal{A}{\rm t}(\mathcal{E})$ of an $\mathcal{O}_X$-module $\mathcal{E}$. For a given Lie algebroid $\mathcal{L}$, a flat $\mathcal{L}$-connection on $\mathcal{E}$ is equivalent to having an $\mathcal{A}{\rm t}(\mathcal{E})$-Lie algebroid structure on $\mathcal{L}$ (see Definitions \ref{L-connections} and \ref{flat}).
\end{Rem}

We now introduce a rooted version of sheaves of Gerstenhaber algebras (see Section~\ref{sheaves of G-algs}).
\begin{Def} Let $\mathcal{A}:=\bigl(\bigoplus^{\infty}_{i=0} \mathcal{A}^i, \wedge, [\cdot, \cdot]\bigr)$ be a Gerstenhaber algebra over $(X, \mathcal{O}_X)$. A Gerstenhaber algebra $\mathcal{A'}:=\bigl(\bigoplus^{\infty}_{i=0} \mathcal{A'}^i, \wedge', [\cdot, \cdot]'\bigr)$ over $(X, \mathcal{O}_X)$ is said to be an \emph{$\mathcal{A}$-Gerstenhaber algebra} if there exists a (sheaf of) Gerstenhaber algebra homomorphism (see Definition~\ref{sheaf of G-algs and morphisms})
\[
\psi\colon \ \mathcal{A'} \rightarrow \mathcal{A}.
\]
\end{Def}
	Let $\mathcal{E}$ be an $\mathcal{O}_X$-module. Consider the exterior algebra $\wedge_{\mathcal{O}_X}^\bullet \mathcal{E}:=\bigl(\bigoplus^{\infty}_{i=0} \wedge_{\mathcal{O}_X}^i \mathcal{E}, \wedge\bigr)$, a graded $\mathcal{O}_X$-algebra. We obtain the following results.

\begin{Exm}\quad
	\begin{enumerate}\itemsep=0pt
\item Let $(\mathcal{L'},\phi)$ be an $\mathcal{L}$-Lie algebroid over $(X,\mathcal{O}_X)$. Then the induced map
\[
\wedge_{\mathcal{O}_X}^\bullet \phi\colon \ \mathcal{G}_{\mathcal{L'}}:=\bigl(\wedge_{\mathcal{O}_X}^\bullet \mathcal{L'}, [\cdot,\cdot]^{\mathcal{L'}}\bigr) \rightarrow \mathcal{G}_{\mathcal{L}}:=\bigl(\wedge_{\mathcal{O}_X}^\bullet \mathcal{L}, [\cdot,\cdot]^{\mathcal{L}}\bigr)
\]
 provides a $\mathcal{G}_{\mathcal{L}}$-Gerstenhaber algebra structure on $\mathcal{G}_{\mathcal{L'}}$ (see Example~\ref{Lie alg to G-alg}).
		\item A Gerstenhaber algebra structure on $\wedge^\bullet_{\mathcal{O}_X} \mathcal{E}$ can be viewed as a $\mathcal{G}_{\mathcal{T}_X}$-Gerstenhaber algebra.
	\end{enumerate}
\end{Exm}

\begin{Rem}
	Let $\mathcal{A}':=\bigl(\bigoplus^{\infty}_{i=0} \mathcal{A}'^i, \wedge', [\cdot, \cdot]'\bigr)$ be an $\mathcal{A}:=\bigl(\bigoplus^{\infty}_{i=0} \mathcal{A}^i, \wedge, [\cdot, \cdot]\bigr)$-Gerstenhaber algebra over $(X, \mathcal{O}_X)$. Then by the analogy of Remark~\ref{G-alg to Lie alg}, we get a canonical $\mathcal{A}^1$-Lie algebroid structure on~$\mathcal{A}'^1$ over $(X,\mathcal{O}_X)$.
\end{Rem}
In the following, we extend Theorem~\ref{G-alg & Lie alg} to the framework of $\mathcal{L}$-Lie algebroids.
\begin{Thm} \label{L-version of G-alg & Lie alg}
For an $\mathcal{O}_X$-module $\mathcal{E}$ and a Lie algebroid $\mathcal{L}$, the $\mathcal{O}_X$-algebra $\wedge_{\mathcal{O}_X}^\bullet \mathcal{E}$ has a $\mathcal{G}_{\mathcal{L}}$-Gerstenhaber algebra structure if and only if the $\mathcal{O}_X$-module $\mathcal{E}$ has an $\mathcal{L}$-Lie algebroid structure.
\end{Thm}

\begin{proof}
	Let the $\mathcal{O}_X$-algebra $\wedge_{\mathcal{O}_X}^\bullet \mathcal{E}$ has a $\mathcal{G}_{\mathcal{L}}$-Gerstenhaber algebra structure.
	Thus, $\wedge_{\mathcal{O}_X}^\bullet \mathcal{E}$ is equipped with a Gerstenhaber bracket
\[
[\cdot,\cdot]^{\mathcal{E}}\colon \ \wedge_{\mathcal{O}_X}^i \mathcal{E} \otimes \wedge_{\mathcal{O}_X}^j \mathcal{E} \rightarrow \wedge_{\mathcal{O}_X}^{i+j-1} \mathcal{E},
\]
	for all $i,j \geq 0$ (note that $\wedge_{\mathcal{O}_X}^{-1} \mathcal{E}:=\{0\}$ and $\wedge_{\mathcal{O}_X}^{0} \mathcal{E}:=\mathcal{O}_X$) and there is a Gerstenhaber algebra homomorphism
\[
\psi\colon \ \left(\bigoplus^{\infty}_{i=0} \wedge_{\mathcal{O}_X}^i \mathcal{E}, \wedge, [\cdot,\cdot]^{\mathcal{E}}\right) \rightarrow \left(\bigoplus^{\infty}_{i=0} \wedge_{\mathcal{O}_X}^i \mathcal{L}, \wedge, [\cdot,\cdot]^{\mathcal{L}}\right).
\]
	Our claim is the triple $\bigl(\mathcal{E}, [\cdot,\cdot]_{\mathcal{E}}:= [\cdot,\cdot]^{\mathcal{E}}|_{\mathcal{E}}, \mathfrak{a'}:=\mathfrak{a}\circ \psi|_{\mathcal{E}}\bigr)$ forms an $\mathcal{L}$-Lie algebroid. Since the homomorphism ${\psi}$ is a degree zero map (or, a graded map) and thus the restriction map
\[
\phi:=\psi|_{\mathcal{E}}\colon \ (\mathcal{E}, [\cdot,\cdot]_{\mathcal{E}}) \rightarrow (\mathcal{L}, [\cdot,\cdot]_{\mathcal{L}})
\]
 is an $\mathcal{O}_X$-module homomorphism and $\mathbb{K}_X$-Lie algebra homomorphism. For any two section $D,D' \in \mathcal{E}$ and a section $f \in \mathcal{O}_X$, we obtain the Leibniz identity
\[
[D, f D']_{\mathcal{E}}=[D, f \wedge D']^{\mathcal{E}}
		=[D,f]^{\mathcal{E}}\wedge D'+ f\wedge [D,D']^{\mathcal{E}}
		=\mathfrak{a'}(D)(f) D'+f [D,D']_{\mathcal{E}}
\]
	(since ${\mathfrak{a}'}(D)(f)=\mathfrak{a}(\phi(D))(f)=[\phi(D),f]^{\mathcal{L}}=[\psi(D),\psi(f))]^{\mathcal{L}} =\psi\bigl([D,f]^{\mathcal{E}}\bigr)=[D,f]^{\mathcal{E}}$ holds).
	
	Let $(\mathcal{E}, \mathfrak{a'},[\cdot,\cdot]_{\mathcal{E}})$ be an $\mathcal{L}$-Lie algebroid. Thus, we have a Lie algebroid homomorphism
\[
\phi\colon \ (\mathcal{E}, \mathfrak{a'},[\cdot,\cdot]_{\mathcal{E}}) \rightarrow (\mathcal{L}, \mathfrak{a},[\cdot,\cdot]_{\mathcal{L}}).
\]
	We show that the homomorphism $\phi$ lifts or extends (uniquely) to a Gerstenhaber algebra homomorphism
\[
\tilde{\phi}\colon \ \left(\bigoplus^{\infty}_{i=0} \wedge_{\mathcal{O}_X}^i \mathcal{E}, \wedge, [\cdot,\cdot]^{\mathcal{E}}\right) \rightarrow \left(\bigoplus^{\infty}_{i=0} \wedge_{\mathcal{O}_X}^i \mathcal{L}, \wedge, [\cdot,\cdot]^{\mathcal{L}}\right)
\]
	between the canonical Gerstenhaber algebras associated with the Lie algebroids $\mathcal{E}$ and $\mathcal{L}$ respectively (see Example~\ref{Lie alg to G-alg}).
	Consider two homogeneous section $\tilde{D}, \tilde{D'} \in \wedge_{\mathcal{O}_X}^\bullet \mathcal{E}$, i.e., there is $f,g \in \mathcal{O}_X$ and $D_1, \dots, D_k, D'_1, \dots, D'_l \in \mathcal{E}$ such that
	$\tilde{D}=f D_1\wedge \cdots \wedge D_k$, $\tilde{D'}=g D'_1\wedge \cdots \wedge D'_l$.
	Therefore, the map (sheaf homomorphism over $(X, \mathcal{O}_X)$) $\tilde{\phi}$ is defined as follows.
\begin{gather*}
\tilde{\phi}\bigl(\tilde{D} \wedge \tilde{D'}\bigr):=f \phi(D_1)\wedge \cdots \wedge \phi(D_k)\wedge g \phi(D'_1)\wedge \cdots \wedge \phi(D'_l),\\
\tilde{\phi}\bigl([\tilde{D},\tilde{D'}]^{\mathcal{E}}\bigr):=[f \phi(D_1)\wedge \cdots \wedge \phi(D_k), g \phi(D'_1)\wedge \cdots \wedge \phi(D'_l)]^{\mathcal{L}}.
\end{gather*}
	Hence, the map $\tilde{\phi}$ is a graded $\mathcal{O}_X$-algebra and graded $\mathbb{K}_X$-Lie algebra homomorphism.
\end{proof}
\begin{Rem}
	Recall that the canonical correspondence between the category of Lie--Rinehart algebras and the category of Gerstenhaber algebras is established via standard adjoint functors (see \cite{MG-SS, JH2}). Its algebro-geometric counterpart is illustrated in Example~\ref{Lie alg to G-alg} and Remark~\ref{G-alg to Lie alg}. Consequently, Theorem~\ref{L-version of G-alg & Lie alg} can also be interpreted within this categorical framework.
\end{Rem}

\section[Omega\^{}bullet\_L-differential graded algebras]{$\boldsymbol{\Omega^{\bullet}_{\mathcal{L}}}$-differential graded algebras}	\label{Sec 4}

Lie algebroids give rise to differential graded algebras, and conversely, a specific type of differential graded algebras can define Lie algebroids (see Theorem~\ref{Lie alg & dga}). To generalize this canonical correspondence in the context of $\mathcal{L}$-Lie algebroids, we introduce the concept of $\Omega^{\bullet}_{\mathcal{L}}$-differential graded algebras. This notion extends the framework of (sheaves of) differential graded algebras, enabling the study of $\mathcal{L}$-Lie algebroid structures on locally free $\mathcal{O}_X$-modules of finite rank, where $(X, \mathcal{O}_X)$ is a special ringed space (see Remark~\ref{special spaces}).
\begin{Rem}
	For a locally free Lie algebroid $\mathcal{L}$ of finite rank over $(X,\mathcal{O}_X)$, we have the standard differential graded algebra $\Omega^{\bullet}_{\mathcal{L}}:=\bigl(\wedge^{\bullet}_{\mathcal{O}_X}\mathcal{L}^*,{\rm d}_{\mathcal{L}}\bigr)$ of the Chevalley--Eilenberg--de Rham complex $($see Section~\ref{C-E-d cohomology}$)$.
\end{Rem}
\begin{Def} Let $\mathcal{A}^{\bullet}:=\bigl(\bigoplus^{\infty}_{i=0} \mathcal{A}^i, \wedge, {\rm d}_{\mathcal{A}}\bigr)$ be a differential graded algebra over $(X,\mathcal{O}_X)$.
	It is said to be a~\emph{$\Omega^{\bullet}_{\mathcal{L}}$-differential graded algebra} if for a given Lie algebroid $\mathcal{L}$ over $(X,\mathcal{O}_X)$, there is a differential graded algebra homomorphism
\[
\psi\colon \ \Omega^{\bullet}_{\mathcal{L}} \rightarrow \mathcal{A}^{\bullet}
\]
	such that the differential ${\rm d}_{\mathcal{A}}$ is given by the composition $\psi_1 \circ \mathfrak{a}_{\mathcal{L}}^* \circ {\rm d}_X$ of the maps
\[
\mathcal{O}_X \overset{{\rm d}_X}{\longrightarrow}\Omega^1_X \overset{\mathfrak{a}_{\mathcal{L}}^*}{\longrightarrow} \Omega^1_{\mathcal{L}} \overset{\psi_1}{\longrightarrow} \mathcal{A}^1,
\]
	where ${\rm d}_X$ is the differential, $\mathfrak{a}_{\mathcal{L}}^*$ is the dual of the anchor map $\mathfrak{a}_{\mathcal{L}}$ and
	$\psi_1$ is the restriction of $\psi$ on $ \Omega^1_{\mathcal{L}}:= \mathcal{L}^*$.
\end{Def}

If the exterior algebra $\wedge_{\mathcal{O}_X}^\bullet \mathcal{E^*}$ is equipped with a differential ${\rm d}_{\mathcal{E}}\colon \wedge_{\mathcal{O}_X}^\bullet \mathcal{E^*} \rightarrow \wedge_{\mathcal{O}_X}^\bullet \mathcal{E}^{*+1}$, for an $\mathcal{O}_X$-module $\mathcal{E}$, then the following result holds, extends Theorem~\ref{Lie alg & dga}.

\begin{Thm} \label{L-Lie alg & dga}
	Let $(\mathcal{L}, [\cdot, \cdot]_{\mathcal{L}}, \mathfrak{a}_{\mathcal{L}})$ $($or simply $\mathcal{L})$ be a locally free Lie algebroid of finite rank over $(X,\mathcal{O}_X)$.
	A~locally free $\mathcal{O}_X$-module $\mathcal{E}$ of finite rank has an $\mathcal{L}$-Lie algebroid structure if and only if $\Omega^{\bullet}_{\mathcal{E}}:=(\wedge^{\bullet}_{\mathcal{O}_X}\mathcal{E}^*,{\rm d}_{\mathcal{E}})$ has a $\Omega^{\bullet}_{\mathcal{L}}$-differential graded algebra structure. $($For more generalities, see Remark~{\rm \ref{dga-reflexive})}.
\end{Thm}
\begin{proof}
	Let the $\mathcal{O}_X$-module $\mathcal{E}$ has an $\mathcal{L}$-Lie algebroid structure given by a Lie algebroid homomorphism
\[
\phi\colon \ (\mathcal{E}, [\cdot, \cdot]_{\mathcal{E}}, \mathfrak{a}_{\mathcal{E}}) \rightarrow (\mathcal{L}, [\cdot, \cdot]_{\mathcal{L}}, \mathfrak{a}_{\mathcal{L}}).
\]
	We consider the Chevalley--Eilenberg--de Rham complex $\Omega^{\bullet}_{\mathcal{E}}:=\bigl(\wedge^{\bullet}_{\mathcal{O}_X}\mathcal{E}^*,{\rm d}_{\mathcal{E}}\bigr)$ associated to the Lie algebroid $(\mathcal{E}, [\cdot, \cdot]_{\mathcal{E}}, \mathfrak{a}_{\mathcal{E}})$ $($or simply $\mathcal{E})$. Note that, the anchor $\mathfrak{a}_{\mathcal{E}}=\mathfrak{a}_{\mathcal{L}} \circ \phi$ and the differential~${\rm d}_{\mathcal{E}}$ is given by the composition of the maps
\[
\mathcal{O}_X \overset{{\rm d}_X}{\longrightarrow}\Omega^1_X \overset{\mathfrak{a}_{\mathcal{L}}^*}{\longrightarrow} \Omega^1_{\mathcal{L}}:= \mathcal{L}^* \overset{\phi^*}{\longrightarrow} \mathcal{E}^*=:\Omega^1_{\mathcal{E}},
\]
 i.e., ${\rm d}_{\mathcal{E}}$ is the canonical extension of the map $\phi^* \circ \mathfrak{a}_{\mathcal{L}}^* \circ {\rm d}_X$, where
\[
{\rm d}_{\mathcal{E}}f (D)={\rm d}_Xf(\mathfrak{a}_{\mathcal{L}}(\phi(D)))= \mathfrak{a}_{\mathcal{L}}(\phi(D))(f)\qquad \text{for all $D \in \mathcal{E}$ and $f \in \mathcal{O}_X$}.
\]
We show that the induced map forms a homomorphism of the (sheaves of) cochain complexes
\[
\widetilde{\phi^*}\colon \ \Omega^{\bullet}_{\mathcal{L}} \rightarrow \Omega^{\bullet}_{\mathcal{E}}
\] is defined by
	$\widetilde{\phi^*}(\omega_1 \wedge \cdots \wedge \omega_k)= (\omega_1 \circ \phi )\wedge \cdots \wedge (\omega_k \circ \phi )$,
	for $\omega_1, \dots, \omega_k \in \mathcal{L}^*$ as follows. Consider $\tilde{\phi}:=\wedge^k \phi$, thus
	for a section $\omega$ of $ \Omega^{k}_{\mathcal{L}}:= \wedge^k_{\mathcal{O}_X} \Omega^{1}_{\mathcal{L}}$ and sections $s_1, \dots, s_{k+1}$ of $\mathcal{E}$ we get the following identities	(see (\ref{C-E-D differential}))	
	\begin{align*}
			& 	{\rm d}_{\mathcal{E}}\bigl(\widetilde{\phi^*}(\omega)\bigr)(s_1 \wedge \cdots \wedge s_{k+1})\\
			&\qquad = {\rm d}_{\mathcal{E}}\bigl(\omega \circ \tilde{\phi}\bigr)(s_1 \wedge \cdots \wedge s_{k+1})\\
			&\qquad = \sum^{k+1}_{i=1}(-1)^{i+1}\mathfrak{a}_{\mathcal{E}}(s_i)\bigl(\bigl(\omega \circ \tilde{\phi}\bigr) (s_1\wedge \cdots \wedge \widehat{s_i} \wedge \cdots \wedge s_{k+1})\bigr)\\
			&\qquad\quad
			+ \sum_{i< j}(-1)^{i+j} \bigl(\omega \circ \tilde{\phi}\bigr)([s_i, s_j]_{\mathcal{E}} \wedge s_1\wedge \cdots \wedge \widehat{s_i}\wedge \cdots \wedge \widehat{s_j}\wedge \cdots \wedge s_{k+1})\\
			&\qquad = \sum^{k+1}_{i=1}(-1)^{i+1}\mathfrak{a}_{\mathcal{L}}(\phi(s_i))\bigl(\omega \bigl(\phi(s_1) \wedge \cdots \wedge \widehat{\phi(s_i)} \wedge \cdots \wedge \phi(s_{k+1})\bigr)\bigr)\\
			&
			\qquad\quad + \sum_{i< j}(-1)^{i+j} \omega \bigl([\phi(s_i), \phi(s_j)]_{\mathcal{L}} \wedge \phi(s_1) \wedge \cdots \wedge \widehat{\phi(s_i)} \wedge \cdots \wedge \widehat{\phi(s_j)}\wedge \cdots \wedge \phi(s_{k+1})\bigr)\\
			&\qquad= {\rm d}_{\mathcal{L}}(\omega) (\phi(s_1) \wedge \cdots \wedge \phi(s_{k+1}))
			 = \widetilde{\phi^*}({\rm d}_{\mathcal{L}}(\omega)) (s_1 \wedge \cdots \wedge s_{k+1})
	\end{align*}
	(this is an extension of the identity $\phi^*({\rm d}_{\mathcal{L}}(f))(s)={\rm d}_{\mathcal{L}}(f)(\phi(s))$ where $f \in \mathcal{O}_X$ and $s \in \mathcal{E}$),
	i.e., the required compatibility condition ${\rm d}_{\mathcal{E}} \circ \widetilde{\phi^*}= \widetilde{\phi^*} \circ {\rm d}_{\mathcal{L}}$ holds. By the construction, the map $\widetilde{\phi^*}$ is also a graded algebra map.
	Thus, it provides a $\Omega^{\bullet}_{\mathcal{L}}$-differential graded algebra structure on $\Omega^{\bullet}_{\mathcal{E}}$.
	
Suppose, $\Omega^{\bullet}_{\mathcal{E}}:=\bigl(\wedge^{\bullet}_{\mathcal{O}_X}\mathcal{E}^*,{\rm d}_{\mathcal{E}}\bigr)$ is a $\Omega^{\bullet}_{\mathcal{L}}$-differential graded algebra with connecting morphism $\psi\colon \Omega^{\bullet}_{\mathcal{L}} \rightarrow \Omega^{\bullet}_{\mathcal{E}}$.
	Consider the restriction map $ \psi|_{\mathcal{L}^*}\colon \mathcal{L}^* \rightarrow \mathcal{E}^*$ and dualizing it we get a homomorphism $\phi\colon \mathcal{E}\rightarrow \mathcal{L}$ of $\mathcal{O}_X$-modules, since both $\mathcal{E}$ and $\mathcal{L}$ are locally free $\mathcal{O}_X$-modules of finite rank (see Remark~\ref{loc free f.r.}).
	Define the $\mathcal{O}_X$-module homomorphism
	$\mathfrak{a}_{\mathcal{E}}\colon \mathcal{E}\rightarrow \mathcal{T}_X$ as
\[
\mathfrak{a}_{\mathcal{E}}(D)(f):={\rm d}_{\mathcal{E}}f (D)
\]
	(well definedness of the map is given by the derivation property of the differential ${\rm d}_{\mathcal{E}}$).
	Since ${\rm d}_{\mathcal{E}}=\psi|_{\mathcal{L}^*} \circ \mathfrak{a}_{\mathcal{L}}^* \circ {\rm d}_X= (\mathfrak{a}_{\mathcal{L}}\circ \phi)^* \circ {\rm d}_X$ and thus $\mathfrak{a}_{\mathcal{E}}(D)(f)= (\mathfrak{a}_{\mathcal{L}}\circ \phi)^* ({\rm d}_Xf)(D)={\rm d}_Xf((\mathfrak{a}_{\mathcal{L}}\circ \phi)(D))=((\mathfrak{a}_{\mathcal{L}}\circ \phi)(D))(f)$ for all $f \in \mathcal{O}_X$, $D\in \mathcal{E}$. Hence, the map
	$\mathfrak{a}_{\mathcal{E}}$ should split as $\mathfrak{a}_{\mathcal{L}}\circ \phi$.
	Consider any two sections $D_1, D_2 \in \mathcal{E}$ and a section $\omega \in \mathcal{E}^*$, the Lie bracket $[\cdot,\cdot]_{\mathcal{E}}$ on $\mathcal{E}$ is given by the following identity:
\[
\omega([D_1, D_2]_{\mathcal{E}}):= \mathfrak{a}_{\mathcal{E}}(D_1)(\omega(D_2))- \mathfrak{a}_{\mathcal{E}}(D_2)(\omega(D_1))- {\rm d}_{\mathcal{E}}\omega (D_1 \wedge D_2).
\]
	The map $\mathfrak{a}_{\mathcal{E}}$ is a $\mathbb{K}_X$-Lie algebra homomorphism as follows:
\begin{align*}
\mathfrak{a}_{\mathcal{E}}([D_1,D_2]_{\mathcal{E}})(f)&={\rm d}_{\mathcal{E}}f([D_1,D_2]_{\mathcal{E}}) =\mathfrak{a}_{\mathcal{E}}(D_1)({\rm d}_{\mathcal{E}}f(D_2))-\mathfrak{a}_{\mathcal{E}}(D_2)({\rm d}_{\mathcal{E}}f(D_1))\\
& =\mathfrak{a}_{\mathcal{E}}(D_1)(\mathfrak{a}_{\mathcal{E}}(D_2)(f)) -\mathfrak{a}_{\mathcal{E}}(D_2)(\mathfrak{a}_{\mathcal{E}}(D_1)(f)) =[\mathfrak{a}_{\mathcal{E}}(D_1),\mathfrak{a}_{\mathcal{E}}(D_2)]_c(f)
\end{align*}
	for all $D_1,D_2 \in \mathcal{E}$ and for all $f \in \mathcal{O}_X$.
	To show the Leibniz identity for the map $\mathfrak{a}_{\mathcal{E}}$, we prove that for $f \in \mathcal{O}_X$ and $D_1,D_2 \in \mathcal{E}$ the following identity holds for all $\omega \in \mathcal{E}^*$:
	\begin{align} \label{Leiniz dga}
		\omega([D_1, f D_2]_{\mathcal{E}})= \omega(f [D_1,D_2]_{\mathcal{E}}+\mathfrak{a}_{\mathcal{E}}(D_1)(f) D_2).
	\end{align}
	Notice that, $(\mathcal{E}, [\cdot, \cdot]_{\mathcal{E}}, \mathfrak{a}_{\mathcal{E}})$ forms an $\mathcal{L}$-Lie algebroid, using the following identities (implies identity~(\ref{Leiniz dga})):
	\begin{align*}
			\omega([D_1, f D_2]_{\mathcal{E}})
			&= \mathfrak{a}_{\mathcal{E}}(D_1)(f \omega(D_2))- f \mathfrak{a}_{\mathcal{E}}(D_2)(\omega(D_1))-f {\rm d}_{\mathcal{E}}\omega (D_1 \wedge D_2)\\
			&=
			\mathfrak{a}_{\mathcal{E}}(D_1)(f) \omega(D_2)\hspace{-0.8pt}+\hspace{-0.8pt} f \mathfrak{a}_{\mathcal{E}}(D_1)(\omega(D_2)))\hspace{-0.8pt}-\hspace{-0.8pt} f \mathfrak{a}_{\mathcal{E}}(D_2)(\omega(D_1))\hspace{-0.8pt}-\hspace{-0.8pt}f {\rm d}_{\mathcal{E}}\omega (D_1 \hspace{-0.5pt}\wedge\hspace{-0.5pt} D_2)\\
			&=
			f (\mathfrak{a}_{\mathcal{E}}(D_1)(\omega(D_2))- \mathfrak{a}_{\mathcal{E}}(D_2)(\omega(D_1))-{\rm d}_{\mathcal{E}}\omega (D_1 \wedge D_2))+\mathfrak{a}_{\mathcal{E}}(D_1)(f) \omega(D_2)\\
			&=
			f \omega([D_1, D_2]_{\mathcal{E}}) + \mathfrak{a}_{\mathcal{E}}(D_1)(f) \omega(D_2)
			=\omega(f [D_1,D_2]_{\mathcal{E}}+\mathfrak{a}_{\mathcal{E}}(D_1)(f) D_2).\tag*{\qed}
	\end{align*}\renewcommand{\qed}{}	
\end{proof}
\begin{Rem} \label{dga-reflexive}
	The condition of locally free $\mathcal{O}_X$-modules of finite rank in the above result can be relaxed by considering both $\mathcal{O}_X$-modules $\mathcal{E}$ and $\mathcal{L}$ as reflexive sheaves. Under this assumption, in the second part of the proof, the restriction of the DGA morphism
\[
\psi\colon \ \bigl(\mathcal{H}{\rm om}_{\mathcal{O}_X}\bigl(\wedge^\bullet_{\mathcal{O}_X}\mathcal{L}, \mathcal{O}_X\bigr), {\rm d}_{\mathcal{L}}\bigr) \rightarrow \bigl(\mathcal{H}{\rm om}_{\mathcal{O}_X}\bigl(\wedge^\bullet_{\mathcal{O}_X}\mathcal{E}, \mathcal{O}_X\bigr), {\rm d}_{\mathcal{E}}\bigr)
\]
 to the grading $1$ induces the same $\mathcal{O}_X$-module homomorphism $\phi\colon \mathcal{E} \rightarrow \mathcal{L}$, leveraging the reflexivity property. For the first part of the proof, we do not even need the reflexive condition.
\end{Rem}
 \begin{Exm}[logarithmic de Rham complex]
	For a divisor $Y$ in a complex manifold (or an algebraic variety) $X$, consider the sheaf of logarithmic differential $k$-forms (see \cite{CMD, PT})
\[
\Omega^k_X(\log Y):=\wedge^k_{\mathcal{O}_X}\Omega^1_X(\log Y),
\]
where $\Omega^1_X(\log Y)$ is the sheaf of meromorphic $1$-forms on $X$ with simple poles along the divisor~$Y$. There is a canonical differential, namely the log de Rham differential
	\begin{align} \label{log de Rham diff}
		{\rm d}^{\log}_{\rm dR}\colon \ \Omega^k_X(\log Y) \rightarrow \Omega^{k+1}_X(\log Y)
	\end{align}
	for $k \in \mathbb{N} \cup \{0\}$, produces the log de Rham complex
	$\Omega^{\bullet}_X(\log Y)$, forms an $\Omega^{\bullet}_X$-DGA for the de Rham complex $\Omega^{\bullet}_X$ of $X$ by the canonical map
\[
\Omega^{\bullet}_X \hookrightarrow \Omega^{\bullet}_X(\log Y).
\]
\end{Exm}

\section[L-Lie bialgebroids and strong differential G\_L-Gerstenhaber algebras]{$\boldsymbol{\mathcal{L}}$-Lie bialgebroids and strong differential\\ $\boldsymbol{\mathcal{G}_{\mathcal{L}}}$-Gerstenhaber algebras} \label{Sec 5}

There is a one to one correspondence between Lie bialgebroids and certain strong differential Gerstenhaber algebras (see Theorem~\ref{SDGA & LBA}). We extend this to a correspondence between two new notions, namely, $\mathcal{L}$-Lie bialgebroids and strong differential $\mathcal{G}_{\mathcal{L}}$-Gerstenhaber algebras.
In the following considerations, we take Lie algebroids over $(X,\mathcal{O}_X)$ that are locally free $\mathcal{O}_X$-modules of finite rank (see Remark~\ref{BV-reflexive} for a broader perspective).
\begin{Def}
	An \emph{$\mathcal{L}$-Lie bialgebroid} is a pair $(\mathcal{L}', \mathcal{L}'^*)$ of $\mathcal{L}$-Lie algebroids over $(X, \mathcal{O}_X)$ that are in duality. The G-bracket $[\cdot, \cdot]^{\mathcal{L}'}=:[\cdot, \cdot]'$ on $\wedge^\bullet_{\mathcal{O}_X} \mathcal{L}'$ satisfies a compatibility condition expressed using the differential ${\rm d}_{\mathcal{L}^{'*}} =:{\rm d}'_*$ on $ \wedge^\bullet_{\mathcal{O}_X} \mathcal{L}' $ $($since $(\mathcal{L'}^{*})^* \cong \mathcal{L'})$, defined by the Lie algebroid structure of $\mathcal{L}'^*$. Specifically, ${\rm d}'_* $ acts as a derivation of the Lie bracket on $\mathcal{L}'$, meaning that for any sections $D_1$, $D_2$ of $\mathcal{L}'$, the following holds:
\[
{\rm d}'_*[D_1,D_2]'=[{\rm d}'_*D_1, D_2]'+[D_1, {\rm d}'_* D_2]'.
\]
\end{Def}
\begin{Rem}A $\mathcal{T}_X$-Lie bialgebroid is same as a Lie bialgebroid over a non-singular space $X$ (see Section~\ref{Lie bialgebroids}).
\end{Rem}
\begin{Exm} \label{log Poisson}
	Let $(X, \pi)$ be a Poisson manifold and $(\mathcal{O}_X, \{\cdot, \cdot \})$ be the associated sheaf of Poisson algebras. Thus, there is a canonical Lie algebroid homomorphism
	$\tilde{\pi}\colon \Omega^1_X \rightarrow \mathcal{T}_X$ induced from the Poisson structure on $X$.
	Consider a closed Poisson analytic subspace $Y:= V(\mathcal{I})$ of $X$, i.e., the zero locus of an ideal sheaf $\mathcal{I}$ of $\mathcal{O}_X$ satisfying the condition $\{\mathcal{I}, \mathcal{O}_X\} \subset \mathcal{I}$.
	Thus, $Y$ is an $\Omega^1_X$-invariant subspace, i.e., $\mathcal{I}{\rm m}(\tilde{\pi})$ preserves $\mathcal{I}$.
	
The sheaf of meromorphic $1$-forms on $X$ with simple poles along $Y$is denoted by $\Omega^1_X(\log Y)$. In particular, when $Y$ is a Poisson hypersurface with singularities, the anchor map $\tilde{\pi}$ induces a~Lie algebroid homomorphism
	\begin{align*}
		\bar{\pi}\colon \ \Omega^1_X(\log Y) \rightarrow \mathcal{T}_X(-\log Y).
	\end{align*}
	Therefore, $\bigl(\Omega^1_X(\log Y), \mathcal{T}_X(-\log Y)\bigr)$ forms a $\mathcal{T}_X(-\log Y)$-Lie algebroid if in addition $Y$ is free divisor. Moreover, the dual pair forms a triangular Lie bialgebroid (see~\cite{BP}).
\end{Exm}
\begin{Prop}
	An $\mathcal{L}$-Lie bialgebroid $(\mathcal{L'}, \mathcal{L'}^*)$ is equivalent to the $\mathcal{L}$-Lie bialgebroid $(\mathcal{L'}^*, \mathcal{L'})$.
\end{Prop}
\begin{Rem}Let $(\mathcal{L}, \mathcal{L}^*)$ be a Lie bialgebroid over $(X,\mathcal{O}_X)$. Then the anchor map $\mathfrak{a}'\colon \mathcal{L}^* \rightarrow \mathcal{T}_X$ induces the $\mathcal{O}_X$-module homomorphism $(\mathfrak{a}')^*\colon \Omega^1_X \rightarrow \mathcal{L}$, provides an $\Omega^\bullet_X$-DGA structure on $\wedge^\bullet_{\mathcal{O}_X} \mathcal{L}$ as an extension of $(\mathfrak{a}')^* \circ {\rm d}_X$.
\end{Rem}

\begin{Rem}Let $(\mathcal{L}, \mathcal{L}^*)$ be a Lie bialgebroid over $(X, \mathcal{O}_X)$. Then there is a canonical strong differential Gerstenhaber algebra structure on $\mathcal{G}_{\mathcal{L}}:=\bigl(\bigoplus^{\infty}_{i=0}\wedge_{\mathcal{O}_X}^i \mathcal{L}, \wedge, [\cdot,\cdot]^{\mathcal{L}}\bigr)$ over $(X, \mathcal{O}_X)$ given through the differential ${\rm d}_*\colon \wedge^\bullet_{\mathcal{O}_X}\mathcal{L} \rightarrow \wedge^{\bullet+1}_{\mathcal{O}_X}\mathcal{L}$, defined by the Lie algebroid structure of $\mathcal{L}^*$.
\end{Rem}

We now introduce a rooted version of sheaves of strong differential Gerstenhaber algebras (see Section~\ref{strong diff G-algs}).
\begin{Def}
	A \emph{strong differential $\mathcal{G}_{\mathcal{L}}$-Ger\-sten\-haber algebra} is a $\mathcal{G}_{\mathcal{L}}:=\bigl(\wedge_{\mathcal{O}_X}^\bullet \mathcal{L}, [\cdot,\cdot]^{\mathcal{L}}\bigr)$-Ger\-sten\-haber algebra $\mathcal{A}^\bullet:=\bigl(\bigoplus^{\infty}_{i=0} \mathcal{A}^i, \wedge, [\cdot, \cdot]\bigr)$ equipped with a differential
	$\tilde{{\rm d}}$ on $\mathcal{A}^\bullet$ such that $\bigl(\mathcal{A}^\bullet, \tilde{{\rm d}}\bigr)$ forms a strong differential Gerstenhaber algebra over $(X, \mathcal{O}_X)$. In addition, we need to have a~differential graded algebra $($DGA$)$ homomorphism
\[
\psi\colon \ \Omega^\bullet_{\mathcal{L}}=\bigl(\wedge_{\mathcal{O}_X}^\bullet \mathcal{L}^*, {\rm d}_{\mathcal{L}}\bigr) \rightarrow \bigl(\mathcal{A}^\bullet,\tilde{{\rm d}}\bigr).
\]		
\end{Def}
\begin{Exm} \label{log Poisson strong diff G-algs}
	Let $Y$ be a hypersurface of a complex manifold $X$ defined by the vanishing of a global function. Then, there is a canonical surjective DGA morphism from the logarithmic de Rham complex $\Omega^{\bullet}_X(\log Y):=\bigl(\wedge^{\bullet}_{\mathcal{O}_X}\Omega^1_X(\log Y), {\rm d}^{\log}_{\rm dR}\bigr)$ of $Y$ in $X$ (see (\ref{log de Rham diff})) to the de Rham complex of~$X$. On the other hand, if $X$ is a Poisson manifold and $Y$ is a Poisson hypersurface, then there is an inclusion map of DGA's from the logarithmic Poisson complex $\mathcal{T}^{\bullet}_X(-\log Y):=\bigl(\wedge^{\bullet}_{\mathcal{O}_X}\mathcal{T}_X(-\log Y), {\rm d}^{\log}_{\rm Pois}\bigr)$ to the Poisson complex $\mathcal{T}^{\bullet}_X:=\bigl(\wedge^{\bullet}_{\mathcal{O}_X}\mathcal{T}_X, {\rm d}_{\rm Pois}\bigr)$ of $X$ along $Y$ (see \cite{JH, BP}).
	With the canonical strong differential G-algebra structures over $(X, \mathcal{O}_X)$ (using Examples~\ref{Cotantgent Poisson} and~\ref{log Poisson}), we get the commutative diagram
\[
		\begin{tikzcd}
			{\Omega^{\bullet}_X(\log Y)}
			\arrow[r, two heads]
			\arrow[d, "{\wedge^{\bullet}_{\mathcal{O}_X}\bar{\pi}}"]
			&
			{\Omega^{\bullet}_X}
			\arrow[d, "{\wedge^{\bullet}_{\mathcal{O}_X}\tilde{\pi}}"]\\
			{\mathcal{T}^{\bullet}_X(-\log Y)}
			\arrow[r, hook]
			&
			{\mathcal{T}^{\bullet}_X}.
		\end{tikzcd}
\]
\end{Exm}
\begin{Rem}
	A strong differential $\mathcal{G}_{\mathcal{T}_X}$-Gerstenhaber algebra is same as sheaf of strong differential Gerstenhaber algebra over a special space $(X, \mathcal{O}_X)$.	
\end{Rem}
In the following we consider an analogue of Theorem~\ref{SDGA & LBA} in the $\mathcal{L}$-Lie algebroids settings.
\begin{Thm} \label{L-version of Lie bialg & strong G alg}
Let $\mathcal{L}'$ be a locally free $\mathcal{O}_X$-module of finite rank. The pair $(\mathcal{L}', \mathcal{L}'^*)$ forms an $\mathcal{L}$-Lie bialgebroid over $(X,\mathcal{O}_X)$ if and only if there is a strong differential $\mathcal{G}_{\mathcal{L}}$-Gerstenhaber algebra structure on $\wedge^\bullet_{\mathcal{O}_X} \mathcal{L}'$ over $(X,\mathcal{O}_X)$.
\end{Thm}
\begin{proof}
	Let $(\mathcal{L}', \mathcal{L}'^*)$ be an $\mathcal{L}$-Lie bialgebroids over $(X, \mathcal{O}_X)$. The $\mathcal{L}$-Lie algebroid structure on~$\mathcal{L}'$ is given by a Lie algebroid homomorphism $\phi\colon \mathcal{L}' \rightarrow \mathcal{L}$, provides the canonical Gerstenhaber algebra homomorphism
\[
\tilde{\phi}\colon \ \wedge^\bullet_{\mathcal{O}_X} \mathcal{L'} \rightarrow \wedge^\bullet_{\mathcal{O}_X} \mathcal{L},
\]
	and thus $\bigl(\wedge^\bullet_{\mathcal{O}_X} \mathcal{L'},[\cdot,\cdot]^{\mathcal{L}'}\bigr)$ forms a $\mathcal{G}_{\mathcal{L}}$-Gerstenhaber algebra.
	The Chevalley--Eilenberg--de Rham complex $\Omega^\bullet_{\mathcal{L}'^*}$ of the Lie algebroid $\mathcal{L}'^*$ provides the differential graded algebra $\bigl(\wedge^\bullet_{\mathcal{O}_X} \mathcal{L}', {\rm d}'_*\bigr)$, where the differential
\[
{\rm d}'_*\colon \ \wedge^\bullet_{\mathcal{O}_X} \mathcal{L}' \rightarrow \wedge^{\bullet +1}_{\mathcal{O}_X} \mathcal{L}'
\]
 is a derivation with respect to the standard Gerstenhaber bracket $[\cdot,\cdot]^{\mathcal{L}'}$. Now, consider the $\mathcal{L}$-Lie algebroid structure on $\mathcal{L}'^*$, is given by a Lie algebroid homomorphism $\psi\colon \mathcal{L}'^* \rightarrow \mathcal{L}$, provides the $\mathcal{O}_X$-module homomorphism $\psi^*\colon \mathcal{L}^* \rightarrow \mathcal{L}'$. Thus, we have the following sequence of $\mathcal{O}_X$-modules with $\mathbb{K}_X$-linear maps
\[
\mathcal{O}_X \overset{{\rm d}_X}{\longrightarrow}\Omega^1_X := \mathcal{T}_X^* \overset{\mathfrak{a}_{\mathcal{L}}^*}{\longrightarrow} \Omega^1_{\mathcal{L}}:=\mathcal{L}^* \overset{\psi^*}{\longrightarrow} {\mathcal{L}'}.
\]
	It provides the differential graded algebra homomorphism
\[
\wedge^\bullet \psi^*\colon \ \Omega^\bullet_{\mathcal{L}} \rightarrow \bigl(\wedge^\bullet_{\mathcal{O}_X} \mathcal{L'}, {\rm d}'_*\bigr).
\]
	Therefore, $\bigl(\wedge^\bullet_{\mathcal{O}_X} \mathcal{L'}, [\cdot, \cdot]^{\mathcal{L}'}, {\rm d}'_*\bigr)$ forms a strong differential $\mathcal{G}_{\mathcal{L}}$-Gerstenhaber algebras over $(X, \mathcal{O}_X)$.
	
	Let there is a strong differential $\mathcal{G}_{\mathcal{L}}$-Gerstenhaber algebra structure on
	$\wedge^\bullet_{\mathcal{O}_X} \mathcal{L}'$, i.e., there is a~Gerstenhaber bracket $[\cdot, \cdot]$ and differential $\tilde{{\rm d}}$ on $\wedge^\bullet_{\mathcal{O}_X} \mathcal{L}'$ satisfies the following properties. There exists a (sheaf of) Gerstenhaber algebra homomorphism
\[
\bar{\phi}\colon \ \bigl(\wedge^\bullet_{\mathcal{O}_X} \mathcal{L}', [\cdot, \cdot]\bigr) \rightarrow \bigl(\wedge_{\mathcal{O}_X}^\bullet \mathcal{L}, [\cdot,\cdot]^{\mathcal{L}}\bigr)
\]
	and a (sheaf of) differential graded algebra homomorphism
\[
\psi\colon \ \bigl(\wedge_{\mathcal{O}_X}^\bullet \mathcal{L}^*, {\rm d}_{\mathcal{L}}\bigr) \rightarrow \bigl(\wedge^\bullet_{\mathcal{O}_X} \mathcal{L}',\tilde{{\rm d}}\bigr).
\]
	The homomorphism $\bar{\phi}$ provides an $\mathcal{L}$-Lie algebroid structure on $\mathcal{L}'$ given by the restriction map
\[
\phi:=\bar{\phi}|_{\mathcal{L}'}\colon \ (\mathcal{L}', [\cdot, \cdot]|_{\mathcal{L}'}) \rightarrow (\mathcal{L}, [\cdot, \cdot]_{\mathcal{L}})
\]
	and the homomorphism $\psi$ provides an $\mathcal{L}$-Lie algebroid structure on $\mathcal{L}'^*$, given through dualizing the map
	$\psi_1:= \psi|_{\mathcal{L}^*}\colon \mathcal{L}^* \rightarrow \mathcal{L}'$, i.e., given by the map (sheaf homomorphism)
$\psi_1^*\colon \mathcal{L}'^* \rightarrow \mathcal{L}$.
The compatibility condition for the pair $(\mathcal{L}', \mathcal{L}'^*)$ of $\mathcal{L}$-Lie algebroids to become an $\mathcal{L}$-Lie bialgebroids is given by the derivation property of the differential map $\tilde{d}$ with respect to the graded Lie bracket~$[\cdot, \cdot]$.
\end{proof}
\begin{Rem} \label{BV-reflexive}
	It is important to note that the dual pair of Lie algebroids considered here is not required to consist of vector bundles (or locally free sheaves of finite rank), as in \cite{KM, KM-PX, PX}. However, to properly interpret the duality condition in the definition of a Lie bialgebroid, it is essential for the sheaves involved to be reflexive (see \cite{BP}). As a result, Theorem~\ref{L-version of Lie bialg & strong G alg} can be extended to this more general setting.
\end{Rem}

 \section[Generalized L-connections and L-Lie algebroid homology]{Generalized $\boldsymbol{\mathcal{L}}$-connections and $\boldsymbol{\mathcal{L}}$-Lie algebroid homology}
\label{Sec 6}

We extend Lie algebroid connections and their curvature (see Section~\ref{Atiyah and Lie connections}) to a broader context and develop an analogue of Theorem~\ref{BV-alg & flat conn}. Subsequently, we develop a generalization of Lie algebroid homology~\cite{PX} within this framework. These concepts serve as foundational tools for Section~\ref{Sec 8}.

\subsection[Generalized L-connections]{Generalized $\boldsymbol{\mathcal{L}}$-connections}	

Let $(\mathcal{L}', \phi)$ be an $\mathcal{L}$-Lie algebroid. We consider such $\mathcal{L}'$-connection $\nabla'$ on an $\mathcal{O}_X$-module $\mathcal{E}$ which splits by an $\mathcal{L}$-connection $\nabla$ on $\mathcal{E}$ as follows.
\begin{Def}
	An $\mathcal{L}'$-connection $(\mathcal{E}, \nabla')$ is said to be an $(\mathcal{L}', \mathcal{L})$-connection if there exists an $\mathcal{L}$-connection $(\mathcal{E}, \nabla)$ such that $\nabla'=\nabla \circ \phi$.
\end{Def}
More generally, the above definition applies when considering an $\mathcal{I}{\rm m}(\phi)$-connection $(\mathcal{E}, \nabla)$ instead of an $\mathcal{L}$-connection, as $\mathcal{I}{\rm m}(\phi)$ forms a Lie subalgebroid of $\mathcal{L}$.
\begin{Rem}\quad
	\begin{enumerate}\itemsep=0pt
		\item An $\mathcal{L}$-connection $(\mathcal{E}, \nabla)$ provides the $\mathcal{L}'$-connection $(\mathcal{E}, \nabla \circ \phi)$.
		\item For a Lie algebroid $\mathcal{L}$ over $(X, \mathcal{O}_X)$, an $(\mathcal{L}, \mathcal{L})$-connection on an $\mathcal{O}_X$-module $\mathcal{E}$ is same with a usual $ \mathcal{L}$-connection on $\mathcal{E}$.
	\end{enumerate}	
\end{Rem}

\begin{Def}
	An $(\mathcal{L}', \mathcal{L})$-connection $(\mathcal{E}, \nabla')$ is said to be a flat connection if the curvature of the connection $\nabla'$ is zero, i.e., for any $D_1, D_2 \in \mathcal{L}'$
\[
[\nabla'_{D_1}, \nabla'_{D_2}]_c= \nabla'_{[D_1,D_2]}.
\]
\end{Def}
	\begin{Rem}\quad
	\begin{enumerate}\itemsep=0pt
		\item For a flat $\mathcal{L}$-connection $(\mathcal{E}, \nabla)$, the $\mathcal{L}'$-connection $(\mathcal{E}, \nabla \circ \phi)$ is also flat.
		\item Therefore, for a flat $(\mathcal{L}', \mathcal{L})$-connection $(\mathcal{E}, \nabla')$ we get a Lie algebroid homomorphism
\[
\nabla'\colon \ \mathcal{L}' \rightarrow \mathcal{A}{\rm t}(\mathcal{E}),
\]
where $\mathcal{A}{\rm t}(\mathcal{E})$ is the Atiyah algebroid of the $\mathcal{O}_X$-module $\mathcal{E}$. Thus, the $\mathcal{O}_X$-module homomorphism $\nabla\colon \mathcal{L} \rightarrow \mathcal{E}{\rm nd}_{\mathbb{K}_X}(\mathcal{E})$ satisfying the Leibniz rule, restricted on $\mathcal{I}{\rm m}(\phi)$ is flat.
\end{enumerate}	
\end{Rem}

\subsubsection{Generalized holomorphic connections}
Let $(X, J)$ be a generalized complex manifold \cite{MG-GC, NH}. Then the complexified generalized tangent bundle $\mathbb{T}_{\mathbb{C}}X:=(TX \oplus T^*X) \otimes \mathbb{C}$ can be expressed as the direct sum of $\pm i$-eigen bundles $L^{1,0}$ and $L^{0,1}$ of $J$, i.e.,
$\mathbb{T}_{\mathbb{C}}X = L^{1,0} \oplus L^{0,1}.$
Both the eigen bundles are maximal isotropic Courant involutive subbundles of $\mathbb{T}_{\mathbb{C}}X$, and have canonical Lie algebroid structures. Note that, $L^{0, 1}$ can be identified with the conjugate bundle $\bar{L}^{1, 0}$, as well as the dual bundle $\bigl(L^{1, 0}\bigr)^*$.

A generalized holomorphic (GH) connection $D$ on a smooth complex vector bundle $E$ is given by a (smooth complex) Lie algebroid $L^{1, 0}$-connection on $E$ which is flat \cite{NH}. It is described in the manner of Remark~\ref{L-connection using dual}, as a differential operator
\[
{D}\colon \ \Gamma(E) \rightarrow \Gamma\bigl(E \otimes L^{0,1} \bigr) .
\]
Thus, $D$ is a $\mathbb{C}$-linear map satisfying the canonical Leibniz rule
$D(f s)= \bar{{\rm d}}f s+ f Ds$, where $\bar{{\rm d}}$ results from the composition of homomorphisms
\smash{$C^{\infty}(X) \otimes_{\mathbb{R}} \mathbb{C} \overset{{\rm d} \otimes 1}{\rightarrow} \Omega^1(X) \otimes_{\mathbb{R}} \mathbb{C} \overset{\mathfrak{a}^*}{\rightarrow} \Gamma\bigl({L}^{0,1}\bigr) $}. Here, ${\rm d}$~is~the de Rham differential of $X$ and $\mathfrak{a}$ is the anchor map of $L^{1,0}$. Moreover, $D$ satisfies the curvature zero condition.

Equivalently, the differential operator $D$ on $E$
gives rise to a representation $(\mathcal{E}, \nabla)$ of $\mathcal{L}$ (see \cite{BRT, BP}), i.e.,
\[
\nabla\colon \ \mathcal{L} \rightarrow \mathcal{A}{\rm t}(\mathcal{E}),
\]
where $\nabla$ is the dual of ${D}$, $\mathcal{L}=\Gamma_{X}\bigl(L^{1, 0}\bigr)$ and $\mathcal{E}=\Gamma_X(E)$ are the $\mathcal{O}_X$-modules associated with the sheaf of sections $\Gamma_X(-)$ on the vector bundles $L^{1, 0}$ and $E$ on $X$ respectively, where $\mathcal{O}_X:= C^{\infty}_X \otimes_{\mathbb{R}_X} \mathbb{C}_X$.

In order to study deformations of a GC structure on a manifold $X$, Gualtieri considered an element $\epsilon \in \Gamma\bigl(\wedge^2 L^{0,1}\bigr)= H^0\bigl(X, \wedge^2_{\mathcal{O}_X}\mathcal{L}^*\bigr)$ satisfying the Maurer--Cartan equation~\cite{MG-GC}. The corresponding $i$-eigen bundle is $\mathcal{L}^{\epsilon}:=({\rm Id}+ \epsilon) \mathcal{L}$, with the Lie algebroid homomorphism ${\operatorname{pr}_1\colon \mathcal{L}^{\epsilon} \rightarrow \mathcal{L}}$, given by the projection map. Thus, a GH bundle $(E, D)$ naturally produces flat $(\mathcal{L}^{\epsilon}, \mathcal{L})$\nobreakdash-con\-nec\-tion on $\mathcal{E}:=\Gamma_X(E)$ given by $\nabla^{\epsilon}:=\nabla \circ \operatorname{pr}_1$.
Moreover, under a $B$-field transformation
of symmetries~\cite{MG-GC}, a GH bundle produces such connections.

A pair $(E, D)$ consists with a complex smooth vector bundle $E$ with a GH connection $D$ is known as a GH vector bundle. In particular, a GH bundle over a complex manifold is a co-Higgs bundle, as introduced by Hitchin \cite{NH}.

\subsubsection{Equivariant principal bundles}
In the study of equivariant geometry (see \cite{BKN}), for a holomorphic principal $H$-bundle $E_H$ on a~complex manifold $X$, where $X$ is equipped with a holomorphic action of a connected complex Lie group $G$, a subbundle ${\rm At}^G(E_H)$ of the holomorphic vector bundle ${\rm At}(E_H) \oplus (X \times \operatorname{Lie}(G))$ is constructed.
The pair $\bigl({\rm At}^G(E_H), \beta \circ \alpha\bigr)$ forms a holomorphic Lie algebroid over~$X$, where $\alpha\colon {\rm At}^G(E_H) \rightarrow X \times \operatorname{Lie}(G)$ is given by the composition map
\[
{\rm At}^G(E_H) \hookrightarrow {\rm At}(E_H) \oplus (X \times \operatorname{Lie}(G)) \overset{\operatorname{pr}_2}{\rightarrow} X \times \operatorname{Lie}(G)
\]
and $\beta\colon X \times \operatorname{Lie}(G) \rightarrow TX$ is induced by the action of $G$ on $X$. A $G$-connection on $E_H$ is defined to be a holomorphic splitting of an appropriate short exact sequence, given by a section $\gamma\colon X \times \operatorname{Lie}(G) \rightarrow {\rm At}^G(E_H)$ (i.e., $\alpha \circ \gamma= {\rm id}$).
Consider the Lie algebroid formed by the sheaf of sections
\[
\mathcal{A}{\rm t}^G(E_H):= \Gamma_X\bigl({\rm At}^G(E_H)\bigr) \qquad \text{and}\qquad \mathcal{A}{\rm t}(E_H):= \Gamma_X({\rm At}(E_H))
\]
over $(X, \mathcal{O}_X)$, with the morphism $\mathcal{A}{\rm t}^G(E_H) \overset{\operatorname{pr}_1}{\rightarrow}\mathcal{A}{\rm t}(E_H)$. Consider the sheaf of sections of the transformation Lie algebroid $X \times \operatorname{Lie}(G)$ over~$X$ (see Example~\ref{Lie groups})
\[
\mathcal{L}:= \mathcal{O}_X \otimes_{\mathbb{C}_X} \operatorname{Lie}(G)_X,
\]
together with the connection $\nabla:=\operatorname{pr}_1 \circ \gamma$ on $\mathcal{A}{\rm t}(E_H)$.
When $E_H$ is $G$-equivariant, then it admits a flat $G$-connection $\nabla$.
Now, consider another connected complex Lie group $G'$ together with a Lie group homomorphism $\bar{\phi}\colon G' \rightarrow G$. We have an induced action of $G^{\prime}$ on $E_H$ and a Lie algebroid homomorphism $\phi\colon \mathcal{L}' \rightarrow \mathcal{L}$, where
\[
\mathcal{L}':=\mathcal{O}_X \otimes_{\mathbb{C}_X} \operatorname{Lie}(G')_X.
\]
Then, we obtain a flat $(\mathcal{L}', \mathcal{L})$-connection on the $\mathcal{O}_X$-module corresponding to any vector bundle associated to $E_H$. Analogous notions can be constructed in the context of holomorphic (or smooth or algebraic) vector bundles, and more generally for sheaf theoretic settings (see \cite[Section~3.2]{BK-MP}).
 \begin{Exm} \label{log connection}
	Let $(Y, \mathcal{O}_Y)$ be a closed Poisson analytic hypersurface (or Poisson principal divisor) of a Poisson manifold $(X, \mathcal{O}_X)$. Then the canonical Lie algebroid homomorphism $\tilde{\pi}\colon \Omega^1_X \rightarrow \mathcal{T}_X$ induced by the Poisson bivector field $\pi$, extends to a Lie algebroid homomorphism (see Examples~\ref{Cotantgent Poisson} and~\ref{log Poisson})
\[
\bar{\pi}\colon \ \Omega^1_X(\log Y) \rightarrow \mathcal{T}_X(-\log Y).
\]
	Consider $\mathcal{L}':= \Omega^1_X(\log Y)$, $\mathcal{L}:=\mathcal{T}_X(-\log Y)$. For any flat connection $\nabla$ on an $\mathcal{O}_X$-module $\mathcal{E}$ can be restricted to give a logarithmic connection (see Example~\ref{log Atiyah})
\[
\nabla':=\nabla|_{\mathcal{T}_X(-\log Y)}\colon \ \mathcal{T}_X(-\log Y) \rightarrow \mathcal{A}{\rm t}(\mathcal{E})(-\log Y).
\]
	Thus, we get a flat $(\mathcal{L}', \mathcal{L})$-connection $(\mathcal{E}, \nabla' \circ \bar{\pi})$, described by the following diagram:
 \[
 \begin{tikzcd}
	\Omega^1_X(\log Y) && {\mathcal{T}_X(-\log Y)} \\
	&& {\mathcal{A}{\rm t}(\mathcal{E})(-\log Y).}
	\arrow["{\bar{\pi}}", from=1-1, to=1-3]
	\arrow["{\nabla' \circ \bar{\pi}}"', from=1-1, to=2-3]
	\arrow["{\nabla'}", from=1-3, to=2-3]
\end{tikzcd}\]
\end{Exm}
 We now generalize Theorem~\ref{BV-alg & flat conn} to the following broader setting.
\begin{Thm} \label{L-BV alg}
	Let $\mathcal{L}$ be a locally free $\mathcal{O}_X$-module of rank $n$ $($for some $n \in \mathbb{N})$ and $\bigl(\wedge^{\bullet}_{\mathcal{O}_X}\mathcal{L}', [\cdot,\cdot]'\bigr)$ a Gerstenhaber algebra for some $\mathcal{O}_X$-module $\mathcal{L}'$. Then there is an exact Gerstenhaber algebra $($or BV-algebra$)$ structure on $\wedge^{\bullet}_{\mathcal{O}_X}\mathcal{L}$, say $\bigl(\wedge^{\bullet}_{\mathcal{O}_X}\mathcal{L}, [\cdot,\cdot], \partial\bigr)$, together with a Gerstenhaber algebra homomorphism
\[
\psi\colon \ \bigl(\wedge^{\bullet}_{\mathcal{O}_X}\mathcal{L}', [\cdot,\cdot]'\bigr) \rightarrow \bigl(\wedge^{\bullet}_{\mathcal{O}_X}\mathcal{L}, [\cdot,\cdot]\bigr)
\]
	if and only if there is an $\mathcal{L}$-Lie algebroid structure on $\mathcal{L}'$ with a
	flat $(\mathcal{L}',\mathcal{L})$-connection on the $\mathcal{O}_X$-module $\wedge^{n}_{\mathcal{O}_X}\mathcal{L}$.
\end{Thm}
\begin{proof}
	Suppose $\bigl(\wedge^{\bullet}_{\mathcal{O}_X}\mathcal{L}, [\cdot,\cdot], \partial\bigr)$ is an exact Gerstenhaber algebra over $(X, \mathcal{O}_X)$ together with a Gerstenhaber algebra (sheaf homo)morphism
\[
\psi\colon \ \bigl(\wedge^{\bullet}_{\mathcal{O}_X}\mathcal{L}', [\cdot,\cdot]'\bigr) \rightarrow \bigl(\wedge^{\bullet}_{\mathcal{O}_X}\mathcal{L}, [\cdot,\cdot]\bigr).
\]	
	The Gerstenhaber algebra structure on $\wedge^{\bullet}_{\mathcal{O}_X}\mathcal{L}$ provides a canonical Lie algebroid structure on~$\mathcal{L}$ (Remark~\ref{G-alg to Lie alg}). Moreover, the generating operator $\partial\colon \wedge^{\bullet}_{\mathcal{O}_X}\mathcal{L} \rightarrow \wedge^{\bullet -1}_{\mathcal{O}_X}\mathcal{L}$ provides a flat $\mathcal{L}$\nobreakdash-con\-nec\-tion $\nabla$ on the locally free $\mathcal{O}_X$-module $\wedge^{n}_{\mathcal{O}_X}\mathcal{L}$ of rank $1$, denote it as $\bigl(\wedge^{n}_{\mathcal{O}_X}\mathcal{L}, \nabla\bigr)$. These results are stated in Theorems~\ref{G-alg & Lie alg} and~\ref{BV-alg & flat conn}. By Theorem~\ref{L-version of G-alg & Lie alg}, the Gerstenhaber algebra homomorphism~$\psi$	
	provides an $\mathcal{L}$-Lie algebroid structure on $\mathcal{L}'$ through the connecting map $\psi|_{\mathcal{L}'}$. Therefore, we get the flat $(\mathcal{L}', \mathcal{L})$-connection $\bigl(\wedge^{n}_{\mathcal{O}_X}\mathcal{L}, \nabla \circ \psi|_{\mathcal{L}'}\bigr)$, i.e., we get a Lie algebroid homomorphism by the composition of the following maps
\[
\mathcal{L}' \overset{\psi|_{\mathcal{L}'}}{\longrightarrow} \mathcal{L} \overset{\nabla}{\longrightarrow} \mathcal{A}{\rm t}\bigl(\wedge^{n}_{\mathcal{O}_X}\mathcal{L}\bigr).
\]	
	
Suppose the locally free $\mathcal{O}_X$-module $\mathcal{L}$ is endowed with a Lie algebroid structure and $(\mathcal{L}', \phi)$ be an $\mathcal{L}$-Lie algebroid. Suppose there is a flat $(\mathcal{L}', \mathcal{L})$-connection $\nabla'$ on the $\mathcal{O}_X$-module $\wedge^{n}_{\mathcal{O}_X}\mathcal{L}$. Using Theorem~\ref{BV-alg & flat conn}, we get an exact Gerstenhaber algebra structure on $\wedge^{\bullet}_{\mathcal{O}_X}\mathcal{L}$ and by Theorem~\ref{G-alg & Lie alg} we get a Gerstenhaber algebra structure on $\wedge^{\bullet}_{\mathcal{O}_X}\mathcal{L}'$. Therefore, by Theorem~\ref{L-version of G-alg & Lie alg}, we get the canonical Gerstenhaber algebra homomorphism
\[
\wedge^\bullet \phi\colon \ \bigl(\wedge^{\bullet}_{\mathcal{O}_X}\mathcal{L}', [\cdot,\cdot]'\bigr) \rightarrow \bigl(\wedge^{\bullet}_{\mathcal{O}_X}\mathcal{L}, [\cdot,\cdot]\bigr).\tag*{\qed}
\]\renewcommand{\qed}{}
\end{proof}
 \begin{Rem}
Using Remark~\ref{canonical BV-alg} and Theorem~\ref{L-BV alg}, we can conclude that there is a canonical flat $(\mathcal{L}', \mathcal{L})$-connection on the $\mathcal{O}_X$-module $\wedge^{n}_{\mathcal{O}_X}\mathcal{L}$ and a canonical BV-algebra structure on $\wedge^{\bullet}_{\mathcal{O}_X}\mathcal{L}$.
\end{Rem}

\subsection[L-Lie algebroid homology]{$\boldsymbol{\mathcal{L}}$-Lie algebroid homology}

Let $\mathcal{L}'$ be a Lie subalgebroid of $\mathcal{L}$, where $\mathcal{L}$ be a locally free $\mathcal{O}_X$-module of rank $n$ for some $n \in \mathbb{N}$. Thus, $(\mathcal{L}', \iota)$ is an $\mathcal{L}$-Lie algebroid, where $\iota\colon \mathcal{L}' \rightarrow \mathcal{L}$ is the inclusion map.

Consider a flat $(\mathcal{L}', \mathcal{L})$-connection on the $\mathcal{O}_X$-module $\wedge^{n}_{\mathcal{O}_X}\mathcal{L}$ of rank $1$ $($i.e., sheaf of sections of a line bundle over $X)$.
Let $D$ be its corresponding generating operator, and
\[
\partial= (-1)^{n-\bullet} D\colon \ \wedge^{\bullet}_{\mathcal{O}_X} \mathcal{L} \rightarrow \wedge^{\bullet -1}_{\mathcal{O}_X} \mathcal{L}
\]
$($the importance for choosing the sign of $\partial$ is described in \cite{PX}). Then $\partial^2=0$, and we obtain a~chain complex of sheaves. Thus, the associated sheaves of homologies (see Remark~\ref{hyperhomology})
\[
\mathcal{H}_{\bullet}(\mathcal{L}, \nabla):=\mathcal{K}{\rm er} \partial/{\mathcal{I}{\rm m} \partial}
\] is the sheafification of the presheaf
$U \mapsto H_{\bullet} (\mathcal{L}(U), \nabla(U))$
of the associated homology groups on each open set~$U$ of~$X$, considered as a graded $\mathbb{K}$-vector space. Applying the hyperhomology functor, we get
the graded vector space $\mathbb{H}_{\bullet}(\mathcal{L}, \nabla)$ of hyperhomologies.

The $\mathcal{L}$-Lie algebroid $(\mathcal{L}', \iota)$ is said to be $\partial$-invariant if the restriction map of $\partial$
on $\wedge^{\bullet}_{\mathcal{O}_X} \mathcal{L}'$, say~$\partial'$, is a map as follows
\[
\partial'\colon \ \wedge^{\bullet}_{\mathcal{O}_X} \mathcal{L}' \rightarrow \wedge^{\bullet -1}_{\mathcal{O}_X} \mathcal{L}'.
\]
Therefore, we get a graded vector space homomorphism between the associated homologies
\begin{align} \label{homo on homologies}
&\mathbb{H}_{\bullet}(\mathcal{L}', \nabla') \rightarrow \mathbb{H}_{\bullet}(\mathcal{L}, \nabla), \qquad \text{where}\quad \nabla'= \nabla \circ \iota.
\end{align}
 The map (\ref{homo on homologies}) is an isomorphism if $\mathcal{L}'= \mathcal{L}$ and the connections are homotopic (see~\cite{PX}).

\section[U(O\_X,L)-universal enveloping algebroids]{$\boldsymbol{\mathcal{U}(\mathcal{O}_X, \mathcal{L})}$-universal enveloping algebroids} \label{Sec 7}

The universal enveloping algebra of a Lie--Rinehart algebra serves as a fundamental tool for exploring its homological algebra, analogous to the classical case of Lie algebras (see \cite{ JH, JH1, HOM, GR}). Its algebro-geometric counterpart is the universal enveloping algebroid, a concept that plays a central role in understanding the homological algebra of Lie algebroids (see \cite{UB, BRT, CV,MK, AA, BP}). A key ingredient in this is the Poincar\'e--Birkhoff--Witt (PBW) theorem, which has been established for Lie algebroids in a sheaf-theoretic framework (see \cite{AA}). We further generalize this theory to include the case of $\mathcal{L}$-Lie algebroids, which is useful in Section~\ref{Sec 8}.

We begin by recalling the definition of the universal enveloping algebroid of a Lie algebroid.

	\subsection{Universal enveloping algebroid of a Lie algebroid}
Let $(\mathcal{L},[\cdot,\cdot],\mathfrak{a})$ be a Lie algebroid over $(X, \mathcal{O}_X)$. For each open set $U$ of $X$, let $\mathcal{U}(\mathcal{O}_X(U),\mathcal{L}(U))$ be the universal enveloping algebra of the $(\mathbb{K},\mathcal{O}_X(U))$-Lie--Rinehart algebra $\mathcal{L}(U)$, or equivalently, of the Lie--Rinehart pair $(\mathcal{O}_X(U),\mathcal{L}(U))$ (see \cite{JH,MM, GR}). The sheafification of the canonical presheaf
$U \mapsto \mathcal{U}(\mathcal{O}_X(U),\mathcal{L}(U))$, is known as the universal enveloping algebroid of the Lie algebroid $\mathcal{L}$, and denoted by $\mathcal{U}(\mathcal{O}_X, \mathcal{L})$ (see \cite{UB,DRV,AA,BP,TS} for details).

	From the construction of $\mathcal{U}(\mathcal{O}_X, \mathcal{L})$, it is an associative $\mathbb{K}_X$-algebra and $\mathcal{O}_X$-bimodule. Moreover,
there is a canonical $\mathbb{K}_X$-algebra monomorphism $\iota_X\colon \mathcal{O}_X \hookrightarrow \mathcal{U}(\mathcal{O}_X,\mathcal{L})$ and an $\mathcal{O}_X$-linear map $\iota_{\mathcal{L}}\colon \mathcal{L} \rightarrow \mathcal{U}(\mathcal{O}_X,\mathcal{L})$.
The associative $\mathbb{K}_X$-algebra $\mathcal{U}(\mathcal{O}_X, \mathcal{L})$ is generated by $\mathcal{O}_X$ and $\iota_{\mathcal{L}}(\mathcal{L})$ satisfy the following identities:
\begin{align} \label{identities}
	\bar{D} \bar{D'} - \bar{D'} \bar{D}= \overline{[D, D']}, \qquad \bar{D} f - f \bar{D}= \mathfrak{a}(D)(f),
\end{align}
where $D, D' \in \mathcal{L}$, $f \in \mathcal{O}_X$, and $\bar{D}= \iota_{\mathcal{L}}(D)$ for all $D\in \mathcal{L}$.
Hence, the map $\iota_{\mathcal{L}}$ can also be viewed as a $\mathbb{K}_X$-Lie algebra homomorphism.

The $\mathbb{K}_X$-algebra $\mathcal{U}(\mathcal{O}_X,\mathcal{L})$ has a natural filtration of $\mathcal{O}_X$-bimodules:
\[
\mathcal{O}_X = \mathcal{U}_{(0)}(\mathcal{O}_X,\mathcal{L}) \subset \mathcal{U}_{(1)}(\mathcal{O}_X,\mathcal{L}) \subset \mathcal{U}_{(2)}(\mathcal{O}_X,\mathcal{L}) \subset \cdots,
\]
where $\mathcal{U}_{(n)}(\mathcal{O}_X,\mathcal{L})$ is spanned by the sections of $\mathcal{O}_X$ and the powers $\iota_{\mathcal{L}}(\mathcal{L})^m$ for $m = 1,2, \dots ,n$.
The direct sum of the quotient sheaves associated with each pair of consecutive sheaves appearing in the above filtration on $\mathcal{U}(\mathcal{O}_X,\mathcal{L})$ forms a sheaf of graded algebras denoted by $\operatorname{gr}(\mathcal{U}(\mathcal{O}_X,\mathcal{L}))$. Thus,
\[
\operatorname{gr}(\mathcal{U}(\mathcal{O}_X,\mathcal{L})) = \bigoplus_{n\geq 0} \mathcal{U}_{(n)}(\mathcal{O}_X,\mathcal{L})/\mathcal{U}_{(n-1)}(\mathcal{O}_X,\mathcal{L}),
\]
where $\mathcal{U}_{(-1)}(\mathcal{O}_X,\mathcal{L}):= \{0\}$. It is a commutative associative untial $\mathcal{O}_X$-algebra (using the identities (\ref{identities})).

\begin{Rem}\label{Uni prop}
	The universal enveloping algebroid $\mathcal{U}(\mathcal{O}_X,\mathcal{L})$ of a Lie algebroid $\mathcal{L}$ over $(X, \mathcal{O}_X)$ is characterized by the \emph{universal property}: \label{Universal}	
	Let $\mathcal{A}$ be a sheaf of unital associative $\mathbb{K}_X$-algebra with sheaf homomorphisms $\phi\colon \mathcal{O}_X \rightarrow \mathcal{A}$ of $\mathbb{K}_X$-unital algebras and $\psi\colon (\mathcal{L}, [\cdot,\cdot]) \rightarrow (\mathcal{A}, [\cdot,\cdot]_c)$ of $\mathbb{K}_X$-Lie algebras such that $\phi(f)\psi(D) = \psi(fD)$ and $[\psi(D), \phi(f)]_c = \phi(\mathfrak{a}(D)(f))$ holds for $f \in \mathcal{O}_X$ and $D \in \mathcal{L}$. Then, there exists a unique homomorphism of unital $\mathbb{K}_X$-algebras $\widetilde{\psi}\colon \mathcal{U}(\mathcal{O}_X, \mathcal{L}) \rightarrow \mathcal{A}$ such that $\widetilde{\psi} \circ \iota = \phi$ and $\widetilde{\psi} \circ \iota_{\mathcal{L}} = \psi$.
\end{Rem}
\begin{Exm}
For a non-singular special space $(X,\mathcal{O}_X)$, the universal enveloping algebroid is isomorphic to the sheaf of differential operators $\mathcal{D}_X$ on $X$ (i.e., the sheaf of differential operators over $\mathcal{O}_X$, sometimes denoted as $\mathcal{D}{\rm iff}(\mathcal{O}_X))$ \cite{AA, TS}, i.e.,
\[
\mathcal{U}(\mathcal{O}_X, \mathcal{T}_X) \cong \mathcal{D}{\rm iff}(\mathcal{O}_X)=:\mathcal{D}_X.
\]
\end{Exm}

\begin{Exm} \label{sheaf of log diff operators} Given a principal divisor $Y$ in some complex manifold or smooth algebraic variety~$X$, the sheaf of logarithmic derivations and the sheaf of logarithmic differential operators are denoted as $\mathcal{T}_X(-\log Y)$ and $\mathcal{D}_X(-\log Y)$, respectively~\cite{AA}. In the case of a free divisor~$Y$ in~$X$ (i.e., $\mathcal{T}_X(-\log Y)$ is locally free $\mathcal{O}_X$-module \cite{CMD,AA,BP}), we have (after sheafifying the local description given for the module of logarithmic derivations in~\cite{LM}) the isomorphism
\[
\mathcal{U}(\mathcal{O}_X, \mathcal{T}_X(-\log Y)) \cong \mathcal{D}_X(-\log Y).
\]
\end{Exm}
\begin{Exm}
In Example~\ref{Free Lie algebroid}, we consider the path algebroid $\mathcal{P}_X$ of a smooth manifold or smooth algebraic variety $X$, the free Lie algebroid generated by tangent sheaf $\mathcal{T}_X$. It forms a locally free Lie algebroid over $(X,\mathcal{O}_X)$, possibly of infinite rank. The universal enveloping algebroid $\mathcal{U}(\mathcal{O}_X,\mathcal{P}_X) =:\mathbb{D}_X$ is described as sheaf of non-commutative differential operators on~$X$. A $\mathcal{P}_X$-module is simply a vector bundle with a $\mathcal{T}_X$-connection, studied as a $\mathbb{D}_X$-module, describing systems of linear PDEs on the path space (see~\cite{MK}).
\end{Exm}

In the following we define an algebro-geometric analogue of algebra of differential operators on a module (see \cite{HOM, LM,TS} for a local description).

In general, we consider the sheaf of differential operators on an $\mathcal{O}_X$-module $\mathcal{E}$, denote it by~$\mathcal{D}{\rm iff}(\mathcal{E})$ as follows. It forms a sheaf of almost commutative algebras, define inductively as follows:
\[
\mathcal{D}{\rm iff}(\mathcal{E}):= \bigcup_{n \geq 0} \mathcal{D}{\rm iff}^{\leq n}(\mathcal{E}),
\]
where
$\mathcal{D}{\rm iff}^{0}(\mathcal{E}):= \mathcal{E}{\rm nd}_{\mathcal{O}_X}(\mathcal{E})$,
and for $n \geq 1$ consider the $\mathcal{O}_X$-modules
\[
\mathcal{D}{\rm iff}^{\leq n}(\mathcal{E}):= \bigl\{\tilde{D} \in \mathcal{E}{\rm nd}_{\mathbb{K}_X}(\mathcal{E}) \mid \bigl[\tilde{D}, f\bigr]_c \in \mathcal{D}{\rm iff}^{\leq n-1}(\mathcal{E}),\, f \in \mathcal{O}_X\bigr\}.
\]
Note that, for any $\tilde{D} \in \mathcal{D}{\rm iff}^{\leq n}(\mathcal{E})$,
the map $\bigl[\tilde{D}, f\bigr]_c\colon \mathcal{E} \rightarrow \mathcal{E}$ is given by $s \mapsto \tilde{D}(f s) - f \tilde{D}(s)$, for $s \in \mathcal{E}$.

\subsection{On generalized PBW theorem}
\label{Gen PBW}
 If $\mathcal{L}$ is a locally free Lie algebroid over $(X, \mathcal{O}_X)$, then the graded algebra $\operatorname{gr}(\mathcal{U}(\mathcal{O}_{X}, \mathcal{L}))$ is isomorphic to the symmetric algebra $\mathcal{S}_{\mathcal{O}_{X}}\mathcal{L}$ as graded $\mathcal{O}_X$-algebras (refer as generalized PBW theorem, see \cite{MK,AA, BP}), given by the PBW map $($also known as the symmetrization map$)$
\begin{gather*}
\tilde{\psi}_{\mathcal{L}} \colon \ \mathcal{S}_{\mathcal{O}_X}\mathcal{L} \rightarrow \operatorname{gr}(\mathcal{U}(\mathcal{O}_{X}, \mathcal{L})), \qquad
	 D_1 \otimes \cdots \otimes D_k \mapsto \frac{1}{k!} \sum_{\sigma \in S_k} \bar{D}_{\sigma(1)} \cdots \bar{D}_{\sigma(k)},
\end{gather*}
where $D_1, \dots , D_k$ are sections of $\mathcal{L}$ and $\bar{D}_{i}$ is the associated class of $D_i$ for $i= 1, \dots, k$.

\begin{Rem}\label{embed}
For a locally free Lie algebroid $\mathcal{L}$ over $(X, \mathcal{O}_X)$, the homomorphism (as $\mathbb{K}_X$\nobreakdash-alge\-bras and $\mathcal{O}_X$-modules) $\iota_{\mathcal{L}}\colon \mathcal{L} \rightarrow \mathcal{U}(\mathcal{O}_X,\mathcal{L})$ become an embedding.
	As a consequence (using Remark~\ref{Uni prop}), the category of representations of $\mathcal{L}$ and category of left $\mathcal{U}(\mathcal{O}_X, \mathcal{L})$-modules are equivalent (see \cite{UB, BP, AS}).
\end{Rem}
For an $\mathcal{L}$-Lie algebroid $(\mathcal{L}', \phi)$, we get the canonical map (see \cite{UB, AA})
\[
\bar{\phi}:=\mathcal{U}(\phi)\colon \ \mathcal{U}(\mathcal{O}_X, \mathcal{L}') \rightarrow \mathcal{U}(\mathcal{O}_X, \mathcal{L})
\]
using the functorial property of $ \mathcal{U}(\mathcal{O}_X, -)$. For the universal enveloping algebroids, we have a~canonical $\mathbb{K}_X$-algebra and $\mathcal{O}_X$-bimodule structure, with a standard filtration. Considering their graded quotients, induces a canonical map
\[
\operatorname{gr}\bigl(\bar{\phi}\bigr)\colon \ \operatorname{gr}(\mathcal{U}(\mathcal{O}_X, \mathcal{L}')) \rightarrow \operatorname{gr}(\mathcal{U}(\mathcal{O}_X, \mathcal{L})).
\]
Also, we get a canonical map between the associated sheaf of symmetric algebras, denote it by
\[
\tilde{\phi}\colon \ \mathcal{S}{\rm ym}_{\mathcal{O}_X}(\mathcal{L}') \rightarrow \mathcal{S}{\rm ym}_{\mathcal{O}_X}(\mathcal{L}).
\]

\begin{Thm} \label{L-version PBW}
Let $\mathcal{L}$ be a locally free Lie algebroid over $(X, \mathcal{O}_X)$ and $(\mathcal{L}', \phi)$ a locally free $\mathcal{L}$-Lie algebroid. Then the PBW isomorphisms $($of graded $\mathcal{O}_X$-algebras$)$
\[
\psi_{\mathcal{L}'}\colon \ \operatorname{gr}(\mathcal{U}(\mathcal{O}_X, \mathcal{L}')) \rightarrow \mathcal{S}{\rm ym}_{\mathcal{O}_X} (\mathcal{L}'),\qquad
	\psi_{\mathcal{L}}\colon \ \operatorname{gr}(\mathcal{U}(\mathcal{O}_X, \mathcal{L})) \rightarrow \mathcal{S}{\rm ym}_{\mathcal{O}_X} (\mathcal{L}),
\]
 where $\psi_{\mathcal{L}'}:=\bigl(\tilde{\psi}_{\mathcal{L}'}\bigr)^{-1}$ and $\psi_{\mathcal{L}}:=\bigl(\tilde{\psi}_{\mathcal{L}}\bigr)^{-1}$,
	satisfies the compatibility condition $\tilde{\phi} \circ \psi_{\mathcal{L}'}= \psi_{\mathcal{L}} \circ \operatorname{gr}\bigl(\bar{\phi}\bigr)$.
\end{Thm}
\begin{proof}The proof is straightforward by the following facts.
	We need to use the functorial property of $\mathcal{U}(\mathcal{O}_X, -)$ and of $\mathcal{S}{\rm ym}_{\mathcal{O}_X} (-)$, and then applying the generalized PBW Theorem~\ref{Gen PBW} on $\mathcal{L}'$ and on $\mathcal{L}$.
\end{proof}
	\begin{Cor} For a locally free Lie algebroid $\tilde{\mathcal{L}}$, we have a canonical embedding $\iota_{\tilde{\mathcal{L}}}\colon \tilde{\mathcal{L}} \rightarrow \mathcal{U}\bigl(\mathcal{O}_X, \tilde{\mathcal{L}}\bigr)$.
	For locally free Lie algebroids $\mathcal{L}'$ and $\mathcal{L}$, where $(\mathcal{L}', \phi)$ is an $\mathcal{L}$-Lie algebroid, we have the canonical induced maps
\begin{gather*}
\operatorname{gr}(\mathcal{U}(\phi))\colon \ \operatorname{gr}(\mathcal{U}(\mathcal{O}_X, \mathcal{L}')) \rightarrow \operatorname{gr}(\mathcal{U}(\mathcal{O}_X, \mathcal{L})),
\qquad
\tilde{\phi}\colon \ \mathcal{S}{\rm ym}_{\mathcal{O}_X}(\mathcal{L}') \rightarrow \mathcal{S}{\rm ym}_{\mathcal{O}_X}(\mathcal{L}').
\end{gather*}
	Thus, we have the following commutative diagram of $\mathcal{O}_X$-modules:
	\[
	\begin{tikzcd}
		\mathcal{L}' \arrow{r}{\iota_{\mathcal{L}'}}\arrow[swap]{d}{\phi} & \mathcal{U}(\mathcal{O}_X, \mathcal{L}') \arrow{d}{\mathcal{U}(\phi)} \arrow{r}{{\rm gr}} & \operatorname{gr}(\mathcal{U}(\mathcal{O}_X, \mathcal{L}')) \arrow{d}{\operatorname{gr}(\mathcal{U}(\phi))} \arrow{r}{\psi_{\mathcal{L}'}} & \mathcal{S}{\rm ym}_{\mathcal{O}_X}(\mathcal{L}') \arrow{d}{\tilde{\phi}} \\
		\mathcal{L} \arrow{r}{\iota_{\mathcal{L}}} & \mathcal{U}(\mathcal{O}_X, \mathcal{L}) \arrow{r}{{\rm gr}} & \operatorname{gr}(\mathcal{U}(\mathcal{O}_X, \mathcal{L})) \arrow{r}{\psi_{\mathcal{L}}} &
		\mathcal{S}{\rm ym}_{\mathcal{O}_X}(\mathcal{L}).
	\end{tikzcd}
	\]
	Hence, we obtain the following identity of maps, which expresses the compatibility condition:
\[
\tilde{\phi}\circ \psi_{\mathcal{L}'}\circ \operatorname{gr} \circ \iota_{\mathcal{L}'}= \psi_{\mathcal{L}}\circ \operatorname{gr} \circ \iota_{\mathcal{L}} \circ \phi.
\]
\end{Cor}

\section{Homology and cohomology correspondence} \label{Sec 8}
In this section, we extend the duality correspondence between Lie algebroid homology and cohomology \cite{AW-EL,JH1, JH2, PX} to the $\mathcal{L}$-Lie algebroid setting. The essential notions required for this extension are outlined in Remark~\ref{Hypercohomology}, Sections~\ref{Sec 6} and~\ref{Sec 7}.

First, we prove a general result about Atiyah algebroids and then we apply it in our case.
\begin{Lem} \label{Lemma}
	Let $\theta\colon \mathcal{E}' \rightarrow \mathcal{E}$ be an $\mathcal{O}_X$-module monomorphism. Then it induces a Lie algebroid homomorphism
\[
\tilde{\theta}\colon \ \mathcal{A}{\rm t}(\mathcal{E}) \rightarrow \mathcal{A}{\rm t}(\mathcal{E}').
\]
\end{Lem}
\begin{proof}
	Consider the canonical map induced from $\theta$ in the following:
\[
\bar{\theta}\colon \ \mathcal{E}{\rm nd}_{\mathbb{K}_X}(\mathcal{E}) \rightarrow \mathcal{E}{\rm nd}_{\mathbb{K}_X}(\mathcal{E}')
\]
	defined as $\tilde{D}:=\bar{\theta}(D)$, where
	$\theta\bigl(\tilde{D}(s')\bigr)= D(\theta (s'))$ for $s' \in \mathcal{E}'$, $D \in \mathcal{E}{\rm nd}_{\mathbb{K}_X}(\mathcal{E})$. The above map is well defined due to the condition of injectiveness of the map $\theta$. It forms an $\mathcal{O}_X$-linear map which is also a $\mathbb{K}_X$-Lie algebra homomorphism as follows. To show $\bar{\theta}(f D)= \widetilde{fD}= f \tilde{D}$ for any $f \in \mathcal{O}_X$, $D \in \mathcal{E}{\rm nd}_{\mathbb{K}_X}(\mathcal{E})$, notice that
\[
\theta \bigl(\widetilde{fD}(s')\bigr)=fD(\theta(s'))=f \theta\bigl(\tilde{D}\bigr)(s')= \theta\bigl(f\tilde{D}\bigr)(s')
\]
	holds for all $s' \in \mathcal{E}$. On the other hand, to show $\bar{\theta}([D_1, D_2]_c)=\bigl[\bar{D}_1, \bar{D}_2\bigr]_c$ for any $D_1, D_2 \in \mathcal{E}$, notice that
\begin{gather*}
\theta\bigl(\bar{\theta} ([D_1, D_2]_c)(s')\bigr)= [D_1, D_2]_c(\theta(s'))=D_1(D_2(\theta(s')))- D_2(D_1(\theta(s'))),\\
D_1(D_2(\theta(s')))- D_2(D_1(\theta(s')))= D_1\bigl(\theta(\bar{D}_2(s'))\bigr)- D_2(\theta(\bar{D}_1)(s'))\\
\hphantom{D_1(D_2(\theta(s')))- D_2(D_1(\theta(s')))}{}
=\theta \bigl(\bar{D}_1 \bigl(\bar{D}_2(s')\bigr) - \bar{D}_2\bigl(\bar{D}_1(s')\bigr)\bigr),
\end{gather*}
	holds for all $s' \in \mathcal{E}$, and then use the injectiveness of the map $\theta$.

	We want to show that the restriction of the map $\bar{\theta}$ on the $\mathcal{O}_X$-submodule and $\mathbb{K}_X$-Lie subalgebra $\mathcal{A}{\rm t}(\mathcal{E})$ of $\mathcal{E}{\rm nd}_{\mathbb{K}_X}(\mathcal{E})$ defines a Lie algebroid homomorphism
\[
\tilde{\theta}\colon \ \mathcal{A}{\rm t}(\mathcal{E}) \rightarrow \mathcal{A}{\rm t}(\mathcal{E}').
\]
Let $\sigma$ and $\sigma'$ be the anchor maps for $\mathcal{A}{\rm t}(\mathcal{E})$ and $\mathcal{A}{\rm t}(\mathcal{E}')$, respectively, and $D \in \mathcal{A}{\rm t}(\mathcal{E})$. Given any Lie algebroid homomorphism $\phi\colon \mathcal{A}{\rm t}(\mathcal{E}) \rightarrow \mathcal{A}{\rm t}(\mathcal{E}')$, we have the identity $\sigma= \sigma' \circ \phi$.
	
If $\bar{\theta}$ defines such a homomorphism, then $\sigma'_{\tilde{D}}(f)= \bigl(\sigma' \circ \tilde{\theta}\bigr) (D) (f)=\sigma (D)(f)=\sigma_D(f)$ holds for all $D \in \mathcal{A}{\rm t}(\mathcal{E})$, $f \in \mathcal{O}_X$. Thus, our claim is $\tilde{D}=\tilde{\theta}(D) \in \mathcal{A}{\rm t}(\mathcal{E}')$, i.e., $\tilde{D}$ satisfies the following Leibniz rule:
\[
\tilde{D}(f s')= f \tilde{D}(s')+ \sigma_D(f) s'
\]
	for all $f \in \mathcal{O}_X$, $s' \in \mathcal{E}'$.
	Notice that, the following equalities hold:
\begin{gather*}
\theta\bigl(\tilde{D}(f s')\bigr)= D(\theta (f s'))=D(f \theta(s'))= f D(\theta(s')) + \sigma_D(f) \theta(s'),\\
f D(\theta(s')) + \sigma_D(f) \theta(s')= f \theta \bigl(\tilde{D}(s')\bigr)+ \sigma_D(f) \theta(s')= \theta \bigl(f \tilde{D}(s') + \sigma_D(f) s'\bigr).
\end{gather*}
	The injectivity of the map $\theta$ provides the required condition for the well definedness of the restriction map $\tilde{\theta}$. The well-definedness of this map ensures that it is a Lie algebroid homomorphism by the definition of the map $\bar{\theta}$ and the Atiyah algebroid structure.
\end{proof}

\begin{Con} \label{locally free condition on (L', L)}
	Let $(\mathcal{L}', \phi)$ be an $\mathcal{L}$-Lie algebroid such that both $\mathcal{L}'$ and $\mathcal{L}$ are locally free $\mathcal{O}_X$\nobreakdash-module of rank $n$ for some $n \in \mathbb{N}$. Then the following results hold using Lemma~{\rm \ref{Lemma}}.
\end{Con}

 \begin{Rem}
In particular, if we consider $\mathcal{L}'= \mathcal{T}_X(-\log Y)$ and $\mathcal{L}= \mathcal{T}_X$ for some free divisor~$Y$ of a complex manifold (or smooth algebraic variety)~$X$ or if we take $\mathcal{L}'=\Omega^1_X$ and $\mathcal{L}=\mathcal{T}_X$ for a~(smooth or holomorphic) Poisson manifold $X$ (see~\cite{BP}), etc., the above condition holds.
\end{Rem}

\begin{Thm} \label{on homologies}
Let $\nabla$ be a flat $\mathcal{L}$-connection on $\wedge^n_{\mathcal{O}_X} \mathcal{L}$ and the map $\phi\colon \mathcal{L}' \rightarrow \mathcal{L}$ be an injective map.
	Then these induces a flat $(\mathcal{L}', \mathcal{L})$-connection on $\wedge^n_{\mathcal{O}_X} \mathcal{L}'$, provides the graded vector space homomorphism
\[
\mathbb{H}_{\bullet} (\phi)\colon \ \mathbb{H}_{\bullet}(\mathcal{L}', \nabla') \rightarrow \mathbb{H}_{\bullet}(\mathcal{L}, \nabla).
\]
\end{Thm}
\begin{proof} Consider the composition of the following Lie algebroid homomorphism (by Lemma~\ref{Lemma}):
\[
\mathcal{L}' \overset{\phi}{\longrightarrow}\mathcal{L} \overset{\nabla}{\longrightarrow} \mathcal{A}{\rm t}\bigl(\wedge^n_{\mathcal{O}_X} \mathcal{L}\bigr) \overset{\theta}{\longrightarrow} \mathcal{A}{\rm t}\bigl(\wedge^n_{\mathcal{O}_X} \mathcal{L}'\bigr),
\]
where $\theta:= \widetilde{\wedge^n \phi}$ is the induced map associated with the map $\wedge^n \phi\colon \wedge^n_{\mathcal{O}_X} \mathcal{L}' \rightarrow \wedge^n_{\mathcal{O}_X} \mathcal{L}$ on the level of Atiyah algebroids. It provides the flat $(\mathcal{L}', \mathcal{L})$-connection $\nabla':= \theta \circ \nabla \circ \phi$ on $\wedge^n_{\mathcal{O}_X} \mathcal{L}'$.
	
Let the flat $\mathcal{L}$-connection $\nabla$ induces the exact generator $\partial\colon \wedge^{\bullet}_{\mathcal{O}_X} \mathcal{L} \rightarrow \wedge^{\bullet -1}_{\mathcal{O}_X} \mathcal{L}$ and the flat $\mathcal{L}'$-connection $\nabla'$ induces the exact generator $\partial'\colon \wedge^{\bullet}_{\mathcal{O}_X} \mathcal{L}' \rightarrow \wedge^{\bullet -1}_{\mathcal{O}_X} \mathcal{L}'$. We need to show that
\[
\partial \circ \wedge^\bullet \phi = \wedge^{\bullet - 1} \phi \circ \partial'.
\]
	One can verify that $\partial'= \bigl(\wedge^{\bullet - 1} \phi\bigr)^{-1} \circ \partial \circ \wedge^\bullet \phi$ holds, using the isomorphism $\phi^{-1}\colon \mathcal{I}{\rm m}(\mathcal{L}') \rightarrow \mathcal{L}'$.
	
	The required graded vector space homomorphism $\mathbb{H}_{\bullet} (\phi)$ on the level of associated Lie algebroid homologies of $\mathcal{L}'$ and $\mathcal{L}$ with coefficient in $\bigl(\wedge^n_{\mathcal{O}_X} \mathcal{L}', \nabla'\bigr)$ and $\bigl(\wedge^n_{\mathcal{O}_X} \mathcal{L}, \nabla\bigr)$, respectively, follows from the map~$\phi$.
\end{proof}
\begin{Rem}
	Recall that for a flat $\mathcal{L}$-connection $($or an $\mathcal{L}$-module$)$ $(\mathcal{E}, \nabla)$, we define an analogue of the complex (\ref{Chevalley-Eilenberg-de Rham complex}), denoted by $\Omega^{\bullet}_{\mathcal{L}}(\mathcal{E})$, known as the Chevalley--Eilenberg--de Rham complex of $\mathcal{L}$ with coefficient in $\mathcal{E}$, using the connection map $\nabla$ in place of the anchor map $\mathfrak{a}$. The graded algebra of hypercohomologies associated with this complex is denoted by $\mathbb{H}^{\bullet}(\mathcal{L}, \mathcal{E})$ (see \cite{UB, BP, AS}). The standard case described in Section~\ref{C-E-d cohomology} is recovered when the $\mathcal{L}$-module $(\mathcal{O}_X, \mathfrak{a})$ is considered.
\end{Rem}

\begin{Rem} \label{on cohomologies}
	Using Condition \ref{locally free condition on (L', L)}, we have the induced graded vector space homomorphism
\[
\mathbb{H}^{\bullet}(\phi)\colon \ \mathbb{H}^{\bullet}\bigl(\mathcal{L}, \wedge^n_{\mathcal{O}_X} \mathcal{L}\bigr) \rightarrow \mathbb{H}^{\bullet}\bigl(\mathcal{L}', \wedge^n_{\mathcal{O}_X} \mathcal{L}'\bigr).
\]
	The map $\mathbb{H}^{\bullet}(\phi)$ is determined by the homomorphism of cochain complexes
\[
\wedge^\bullet_{\mathcal{O}_X} \phi^* \otimes _{\mathcal{O}_X}\wedge^n_{\mathcal{O}_X} \phi^*\colon \
	\wedge^\bullet_{\mathcal{O}_X} \mathcal{L}^* \otimes_{\mathcal{O}_X}\wedge^n_{\mathcal{O}_X} \mathcal{L} \rightarrow \wedge^\bullet_{\mathcal{O}_X} \mathcal{L}'^* \otimes_{\mathcal{O}_X}\wedge^n_{\mathcal{O}_X} \mathcal{L}',
\]
where the dual map $\phi^*\colon \mathcal{L}^* \rightarrow \mathcal{L}'^*$ of $\phi$ can also be viewed as an $\mathcal{O}_X$-module homomorphism from $\mathcal{L}$ to $\mathcal{L}'$, as both $\mathcal{L}$ and $\mathcal{L}'$ are locally free $\mathcal{O}_X$ modules of finite rank.
\end{Rem}

Rinehart demonstrated that the cohomology of a $(\mathbb{K}, A)$-Lie--Rinehart algebra $L$ with coefficient in an $L$-module $E$, which is projective over~$A$, can be described using the derived functor ${\rm Ext}$ in the category of modules over the universal enveloping algebra $\mathcal{U}(A, L)$ (see \cite{JH, GR}), as
\[
H^{\bullet}(L, E) \cong {\rm Ext}^{\bullet}_{\mathcal{U}(A, L)}(A, E).
\]

Huebschmann employed a dual approach to this concept in an appropriate setting. Consider a~$(\mathbb{K}, A)$-Lie--Rinehart algebra $L$ that is also a projective module of rank~$n$ (for some ${n \in \mathbb{N}}$) over~$A$, together with a BV-algebra structure $\bigl(\wedge^{\bullet}_A L, \partial\bigr)$.
 In this context, Huebschmann introduced flat right $(A, L)$-connections on~$A$, which are equivalent to flat left $(A, L)$-connections~$\nabla$ (induced from the generating operator $\partial$) on $\wedge^n_A L$.
He has defined Lie--Rinehart homology $H_{\bullet}(L, A_{\nabla})$ and expressed it by the derived functor ${\rm Tor}$ in the category of right $\mathcal{U}(A, L)$-modules using the Rinehart resolution of $A$, where $A_{\nabla}$ can be viewed as $\bigl(\wedge^n_A L, \nabla\bigr)$ (see~\cite{JH1}). Explicitly, he showed that
\[
H_{\bullet}(L, A_{\nabla}) \cong {\rm Tor}_{\bullet}^{\mathcal{U}(A, L)}(A_{\nabla}, A).
\]

In the following, we recall the duality theorem for Lie--Rinehart algebras (see \cite[Theorem~4]{JH1}). A very particular case of this theorem was first addressed in the context of smooth Lie algebroids by P.~Xu~\cite{PX}. However, J.~Huebschmann \cite{JH1, JH2} developed a more general formulation for Lie--Rinehart algebras, which does not directly follow from Xu's specific case. To extend these ideas to our setting, we use the description provided for Lie--Rinehart algebras as a local framework.
\begin{Thm} \label{duality L-R}
	For a Lie--Rinehart algebra $(A, L)$ where $L$ is projective $A$-module of rank~$n$ for some $n \in \mathbb{N}$, the homology $H_{\bullet}(L, A_{\nabla})$ of $L$ with coefficient in the flat right $(A, L)$-connec\-tion~$\nabla$ on~$A$, i.e., that of the BV-algebra $\bigl(\wedge^\bullet_{A} {L}, \partial\bigr)$ $($where the exact generator $\partial$ is induced by the flat $L$\nobreakdash-con\-nection $\nabla)$, is naturally isomorphic $($as graded vector spaces$)$ to the cohomology $H^{n-\bullet}\bigl(L, \wedge^n_{A} {L}\bigr)$ of~$L$ with coefficients in the $L$-module $\bigl(\wedge^n_{A} {L}, \nabla\bigr)$.
\end{Thm}

In particular, one can view the case of smooth context, by considering $A=C^\infty(M)$ for some smooth manifold $M$, and using the fact that finitely generated projective $A$-module is equivalent to a smooth vector bundle on $M$. That is, $L$ is the space of global sections $\Gamma\bigl(\tilde{L}\bigr)$ of some Lie algebroid $\tilde{L}$ over $M$. Note that, in \cite{PX}, the duality correspondence described with coefficient taking in the trivial line bundle $\wedge^n {\tilde{L}}$ (or in
free $A$-module $\wedge^n_{A} {L}$ of rank $1$). The case associated with a~nontrivial line bundle as coefficient (twisted coefficients) is not described explicitly, but mentioned the possibility of that generality (see \cite{AW-EL}).

It is well known, assuming certain descriptions of Lie--Rinehart algebras, how to extend these notions to the context of Lie algebroids in algebro-geometric settings (see \cite{UB, CV, MK, JV}). To address homological algebraic aspects, we first sheafifying the descriptions at the level of Lie--Rinehart algebras associated with spaces of sections, and then apply the hyper(co)homology functor (see \cite{UB, BRT, AS}).
Thus, following the approach in \cite{JH1}, we conclude that an analogue of Theorem~\ref{duality L-R} holds for locally free Lie algebroids of finite rank over $(X, \mathcal{O}_X)$, as detailed below. Denote $\mathcal{O}_{X, \nabla}$ by a flat right $\mathcal{L}$-connection $\nabla$ on $\mathcal{O}_X$ (i.e., it is a right $\mathcal{L}$-module).

\begin{Rem} \label{duality for BV-Lie algebroids}
	Let $\mathcal{L}$ be a Lie algebroid of rank $n$, having a BV-algebra structure on $\wedge^{\bullet}_{\mathcal{O}_X} \mathcal{L}$ with exact generator $\partial$. Then the (hyper)homology of $\mathcal{L}$ with coefficient in the right $\mathcal{L}$-module $\mathcal{O}_{X, \nabla}$ (this notion is equivalent to the left $\mathcal{L}$-module \smash{$\bigl(\wedge^{n}_{\mathcal{O}_X} \mathcal{L}, \nabla\bigl)$}) where $\nabla$ induces from $\partial$, is naturally isomorphic to the (hyper)cohomology of $\mathcal{L}$ with coefficient in the $\mathcal{L}$-module \smash{$\bigl(\wedge^{n}_{\mathcal{O}_X} \mathcal{L}, \nabla\bigr)$}. Thus, there is the canonical (duality) isomorphism
\begin{equation}\label{dual iso for Lie}
	\tilde{\psi}_{\mathcal{L}}\colon \
	\mathbb{H}_{\bullet}(\mathcal{L}, \nabla)
	\rightarrow
	\mathbb{H}^{n-\bullet}
 \bigl(\mathcal{L},
	\wedge^n_{\mathcal{O}_X}\mathcal{L} \bigr).
\end{equation}
\end{Rem}

Here, we show an analogue of the duality isomorphism (\ref{dual iso for Lie}) between Lie algebroid homology and Lie algebroid cohomology (see \cite{JH1, JH2, PX}), in the context of $\mathcal{L}$-Lie algebroids.

\begin{Thm} \label{L-version duality} Let $\nabla$ be a flat $\mathcal{L}$-connection on $\wedge^n_{\mathcal{O}_X} \mathcal{L}$ and the map $\phi\colon \mathcal{L}' \rightarrow \mathcal{L}$ be an injective map $($satisfying Condition~{\rm \ref{locally free condition on (L', L)})}.
	Then, the canonical duality isomorphism
\[
\tilde{\psi}_{\mathcal{L}'}\colon \ \mathbb{H}_{\bullet}(\mathcal{L}', \nabla') \rightarrow \mathbb{H}^{n- \bullet}\bigl(\mathcal{L}', \wedge^n_{\mathcal{O}_X} \mathcal{L}'\bigr)
\] factors through the duality isomorphism $\tilde{\psi}_{\mathcal{L}}\colon \mathbb{H}_{\bullet}(\mathcal{L}, \nabla) \cong \mathbb{H}^{n- \bullet}\bigl(\mathcal{L}, \wedge^n_{\mathcal{O}_X} \mathcal{L}\bigr)$ as
\[
\tilde{\psi}_{\mathcal{L}'}=\mathbb{H}^{n- \bullet}(\phi) \circ \tilde{\psi}_{\mathcal{L}} \circ \mathbb{H}_{\bullet} (\phi).
\]
\end{Thm}

\begin{proof}
Let $\mathcal{\bar{L}}$ be a locally free Lie algebroid over $(X, \mathcal{O}_X)$ of rank~$n$. If there is a flat $\mathcal{\bar{L}}$-connection $\bar{\nabla}$ on $\wedge^n_{\mathcal{O}_X} \mathcal{\bar{L}}$, then we get the BV-algebra \smash{$\bigl(\wedge^{\bullet}_{\mathcal{O}_X} \mathcal{\bar{L}}, \bar{\partial}\bigr)$}, where $\bar{\partial}\colon \wedge^{\bullet}_{\mathcal{O}_X} \mathcal{\bar{L}} \rightarrow \wedge^{\bullet -1}_{\mathcal{O}_X} \mathcal{\bar{L}}$ is the BV-operator. An algebro-geometric analogue of the duality isomorphism for finitely generated projective Lie--Rinehart algebra (Theorem~\ref{duality L-R}) is applicable for the Lie algebroid $\mathcal{\bar{L}}$ (Remark~\ref{duality for BV-Lie algebroids}) as follows.
	
Using the locally free finite rank condition, there is an open cover $\{U_x\mid x \in X\}$ of $X$, for which we get a $\mathcal{O}_X|_{U_x}$-module isomorphism $\mathcal{\bar{L}}|_{U_x} \cong \bigoplus^n \mathcal{O}_X|_{U_x}$. Consider the sheafification of the presheaf homomorphism
\[
U \mapsto \bar{\psi}_{\mathcal{\bar{L}}(U)}\colon \ H_{\bullet}\bigl(\mathcal{\bar{L}}(U), \nabla_U\bigr) \rightarrow H^{n-\bullet}\bigl(\mathcal{\bar{L}}(U), \wedge^n_{\mathcal{O}_X(U)} \mathcal{\bar{L}}(U)\bigr),
\]
	which is an isomorphism as we have the local duality isomorphisms $\bar{\psi}_{\mathcal{\bar{L}}(U_x)}$ for each $x \in X$, denote it by $\bar{\psi}_{\bar{\mathcal{L}}}$. That is, the associated (co)homology sheaves are isomorphic as
\[
\bar{\psi}_{\bar{\mathcal{L}}}\colon \ \mathcal{H}_{\bullet}\bigl(\mathcal{\bar{L}}, \nabla\bigr) \overset{\cong }{\rightarrow}\mathcal{H}^{n-\bullet} \bigl(\mathcal{\bar{L}}, \wedge^n_{\mathcal{O}_X} \mathcal{\bar{L}}\bigr).
\]
	Therefore, applying hyperfunctors (see \cite{UB, CW1}) we get the graded vector space isomorphism
\[
\tilde{\psi}_{\bar{\mathcal{L}}}:=\mathbb{H}\bigl(\bar{\psi}_{\mathcal{L}}\bigr)\colon \
 \mathbb{H}_{\bullet}\bigl(\mathcal{\bar{L}}, \nabla\bigr) \overset{\cong }{\rightarrow} \mathbb{H}^{n-\bullet} \bigl(\mathcal{\bar{L}}, \wedge^n_{\mathcal{O}_X} \mathcal{\bar{L}}\bigr).
 \]
Thus, the above duality isomorphism holds for the Lie algebroids $\mathcal{L}$ and $\mathcal{L}'$, denote these by $\tilde{\psi}_{\mathcal{L}}$ and $\tilde{\psi}_{\mathcal{L}'}$, respectively. Using the morphisms described in Theorem~\ref{on homologies} and Remark~\ref{on cohomologies}, we get the required result.
	The following diagram represents the relationship between the isomorphisms~$\tilde{\psi}_{\mathcal{L}'}$ and~$\tilde{\psi}_{\mathcal{L}}$:
\[
	\begin{tikzcd}
		{\mathbb{H}_{\bullet}(\mathcal{L}', \nabla')}
		\arrow[r, "{\tilde{\psi}_{\mathcal{L}'}}"]
		\arrow[d, "{\mathbb{H}_{\bullet}(\phi)}"]&
		{\mathbb{H}^{n-\bullet} \bigl(\mathcal{{L'}}, \wedge^n_{\mathcal{O}_X} \mathcal{{L'}}\bigr)}\\
		{\mathbb{H}_{\bullet}(\mathcal{L}, \nabla)}
		\arrow[r, "{\tilde{\psi}_{\mathcal{L}}}"]&
		{\mathbb{H}^{n-\bullet} \bigl(\mathcal{{L}}, \wedge^n_{\mathcal{O}_X} \mathcal{{L}}\bigr)}
		\arrow[u, "{\mathbb{H}^{n-\bullet}(\phi)}"]
	\end{tikzcd}\tag*{\qed}
\]\renewcommand{\qed}{}
\end{proof}

Such situations can arise from singular foliations, singular subalgebroids, triangular Lie bialgebroids, etc.
\begin{Cor}
If $X$ is a nonsingular or smooth special topological ringed space $($see Remarks~{\rm \ref{special spaces}} and~{\rm \ref{nonsingular})} and $\mathcal{L}' = \mathcal{T}_X = \mathcal{L}$, then the above result recovers an analogue of the Poincar\'e duality isomorphism.
\end{Cor}

\begin{Cor}Lichnerowicz showed that, for a smooth finite-dimensional symplectic manifold, the de Rham and Poisson cohomologies coincide $($see {\rm \cite{JH})}. Thus, for such a manifold $X$, consider the $\mathcal{L}$-Lie algebroid $\mathcal{L}':=\Omega^1_X$, where $\mathcal{L}:=\mathcal{T}_X$. Using Theorem~{\rm \ref{L-version duality}}, we obtain a duality isomorphism at the level of the associated de Rham and Poisson $($co$)$homology hyperfunctors or, hyper$($co$)$homologies $($see Section~{\rm \ref{Hypercohomology})}.
\end{Cor}

\section{Further remarks}
In this article, we have employed certain standard techniques from sheaf theory in obtaining our results.
Here, we briefly recall some of these methods and discuss related possibilities.

\subsection{Local to global settings}

In most of the proofs of our theorems, the underlying arguments rely on techniques from Lie--Rinehart algebras as locally compatible data (producing information on certain presheaves), which are then globalized using the sheafification functor (see Section~\ref{Sec 2}). We mostly work at the level of sheaves of sections, keeping track of the underlying presheaf-level structures. To avoid repetition, we sometimes omit the details of the local-to-global functorial correspondences.

For example, we show that a morphism of Gerstenhaber algebras over $(X,\mathcal{O}_X)$ induces a~morphism of Lie algebroids over $(X,\mathcal{O}_X)$. Equivalently, we obtain a homomorphism of sheaves of Lie--Rinehart algebras from a homomorphism of sheaves of Gerstenhaber algebras.

\emph{Stalkwise description.} All our sheaf-theoretic constructions can be understood in terms of the induced algebraic structures on the stalks. Since the stalks of a presheaf and its associated sheaf coincide, the descriptions involving presheaves or local data are equivalent to those given at the level of stalks (see \cite{UB, SR}).

\emph{\v{C}ech cohomology.} For a nonsingular special space $(X, \mathcal{O}_X)$, one can find good (e.g., affine or Stein) open covers of $X$ that allow standard computations of hypercohomology using \v{C}ech cohomology (see \cite{BRT, BP, AS, PT}). One uses \v{C}ech cohomology to simplify the computations and formulations of the above results.

\subsection{Lie algebroids over schemes}
Our discussions apply to Lie algebroids over Noetherian separated schemes (or schemes of finite type) defined over a field of characteristic zero (see \cite{UB, BRT, MK}).

Let $L$ be a $(\mathbb{K}, A)$-Lie--Rinehart algebra. It gives rise to a Lie algebroid $\mathcal{L}$ over the affine scheme $(X, \mathcal{O}_X)$, where $X = \operatorname{Spec}(A)$ and $\mathcal{O}_X$ denotes the sheaf of regular functions induced by the localization of $A$. By localizing the $A$-module $L$, one obtains the $\mathcal{O}_X$-module~$\mathcal{L}$ (see~\cite{CW}).
Conversely, starting with a Lie algebroid $\mathcal{L}$ over an affine scheme $(X, \mathcal{O}_X)$, we obtain the Lie--Rinehart algebra $\mathcal{L}(X)$ over $\mathcal{O}_X(X)$.
The category of Lie--Rinehart algebras and of quasi-coherent Lie algebroids over an affine scheme are equivalent.

The theory of Lie--Rinehart algebras (see \cite{JH, GR}), or equivalently of Lie algebroids over affine schemes, can be extended to the context of Lie algebroids over certain schemes using the standard local-to-global correspondences (see \cite{UB, BRT, MK, HOM, AA}). Some of these analogies have already been discussed.

Moreover, some of the results can be extended to the broader context of Lie algebroids over ringed sites, which requires adopting a more advanced perspective (see \cite{DRV, CV}).

\subsection{Some applications}

In the following contexts, the rooted or relative versions of Lie algebroids and their associated rooted structures naturally arise and play a crucial role in simplifying several situations.

Let $\mathcal{L}$ be a Lie algebroid over a topological ringed space $(X,\mathcal{O}_X)$.

\emph{Lie algebroid extensions.} To study Lie algebroid extensions of $\mathcal{L}$ by $\mathcal{O}_X$ (or, an $\mathcal{L}$-module) \cite{AS, PT2}, it is natural to work within the subcategory of $\mathcal{L}$-Lie algebroids rather than the whole category of Lie algebroids.

\emph{Lie algebroid connectons.} Likewise, when studying flat $\mathcal{L}$-connections on an $\mathcal{O}_X$-module $\mathcal{E}$, one naturally works with $\mathcal{A}{\rm t}(\mathcal{E})$-Lie algebroid structures on $\mathcal{L}$ (see Remark~\ref{Lie-Atiyah alg}).

\emph{Twisted universal enveloping algebroids.} To study $\mathcal{L}$-connections with non-zero curvature $\omega$, one considers modules over the twisted universal enveloping algebroids $\mathcal{U}(\mathcal{O}_X, \mathcal{L}, \omega)$ (see \cite{HOM, PT2}).
One can view such notions as $\mathcal{U}(\mathcal{O}_X, \mathcal{L})$-universal enveloping algebroids (see Section~\ref{Sec 7}) associated to the $\mathcal{L}$-Lie algebroid $\mathcal{L}_{\omega}$ (see Example~\ref{extensions of L}).

\emph{Deformation theory.} In the study of deformation theory of $\mathcal{L}$-Lie algebroids, one is naturally led to consider $\mathcal{G}_{\mathcal{L}}$-Gerstenhaber algebras and related structures (see \cite{CV, TS}).

\emph{Singular subalgebroids.} To study singular subalgebroids of $\mathcal{L}$ (see \cite{MZ}), sub $\mathcal{G}_{\mathcal{L}}$-Gerstenhaber algebras naturally arise. Following Remark~\ref{G-alg to Lie alg}, one can investigate such subalgebroids via these G-algebras.

R.L.~Klaasse considered several important geometric situations arising from foliation theory, Poisson geometry, symplectic geometry, and generalized complex geometry where such rooted Lie algebroids occur naturally \cite{RK}. These provide a natural setting for studying deformation theory of geometric and algebraic structures via the associated rooted versions of G-algebras and BV-algebras.

Moreover, J.~Huebschmann introduced variants of Lie--Rinehart (LR) algebras, namely twilled LR-algebras, in order to develop an analogue of complex structures (on a $C^\infty$-manifold) in the LR-algebra setting, together with related G-algebras and BV-algebras \cite{twilled-JH}. One can observe that some of these notions also fit naturally into the framework of rooted Lie algebroids.

Therefore, many of the notions introduced and studied in this article may serve as useful ingredients for further developments in these directions.

\subsection*{Acknowledgements}

The authors would like to sincerely thank the editor and the anonymous referees for their careful reading of the manuscript and for their valuable suggestions and comments.
The second author expresses his heartfelt gratitude to Dr.\ Ashis Mandal for introducing him to several foundational concepts that have significantly shaped this work. He also extends his sincere thanks to Dr.\ Satyendra Kumar Mishra for his insightful suggestions. The authors acknowledge the support provided by the Indian Institute of Science Education and Research Pune.

\pdfbookmark[1]{References}{ref}
\LastPageEnding

\end{document}